\documentclass[journal]{IEEEtran}

\usepackage{cite}
\usepackage{amsmath,amssymb,amsfonts,bm}
\usepackage{bbm}
\usepackage{graphicx}
\usepackage{subfigure}
\usepackage{textcomp}
\usepackage{algorithm}
\usepackage{algpseudocode}

\allowdisplaybreaks[4]
\newtheorem{defi}{Definition}
\newtheorem{theo}{Theorem}
\newtheorem{lemma}{Lemma}
\newtheorem{assum}{Assumption}

\newenvironment{proof}{\begin{IEEEproof}}{\end{IEEEproof}}
\begin{document}
	\title{Global Convergence of  Policy Gradient Primal-dual Methods for  Risk-constrained LQRs
		\thanks{Research of the first two authors was supported by National Natural Science Foundation of China under Grant no. 62033006. Research of the third author was supported by the ONR MURI Grant N00014-16-1-2710.}
		\thanks{F. Zhao and K. You are with the Department of Automation and BNRist, Tsinghua University, Beijing 100084, China. e-mail: zhaofr18@mails.tsinghua.edu.cn, youky@tsinghua.edu.cn.}
		\thanks{Tamer~Ba\c{s}ar is with the Coordinated Science Laboratory, University of Illinois at Urbana-Champaign, Urbana, IL 61801 USA. e-mail: basar1@illinois.edu.}
		%\thanks{We would like to thank Kaiqing Zhang for his helpful suggestions.}
		\author{Feiran Zhao, Keyou You, \IEEEmembership{Senior Member,~IEEE}, Tamer~Ba\c{s}ar, \IEEEmembership{Life Fellow,~IEEE}}}
	\maketitle
	
	\begin{abstract} While the techniques in optimal control theory are often model-based, the policy optimization (PO) approach directly optimizes the performance metric of interest. Even though it has been an essential approach for reinforcement learning problems,  there is little theoretical understanding on its performance. In this paper, we focus on the risk-constrained linear quadratic regulator (RC-LQR) problem via the PO approach, which requires addressing a challenging non-convex constrained optimization problem. To solve it, we first build on our earlier result that an optimal policy has a time-invariant affine structure to show that the associated Lagrangian function is coercive, locally gradient dominated and has local Lipschitz continuous gradient, based on which we establish strong duality. Then, we design policy gradient primal-dual methods with global convergence guarantees in both model-based and sample-based settings. Finally, we use samples of system trajectories in simulations to validate our methods.
	\end{abstract}
	
	\begin{IEEEkeywords}
		Risk-constrained LQR, stochastic control, policy optimization, reinforcement learning, gradient descent.
	\end{IEEEkeywords}
	
	\section{Introduction}
	
The techniques in conventional optimal control theory often require an explicit dynamical model. Such a model-based idea is relatively easy to provide theoretical guarantees but is usually sensitive to modeling inaccuracy. Policy optimization (PO) methods, as an end-to-end approach, directly search for an optimal control policy to minimize a performance metric of interest and { has advantages in scenarios where the dynamical model is complex and difficult to identify.} In fact, {it has been proved to be an essential approach for applications of reinforcement learning (RL)~\cite{mnih2015human-level,lillicrap2016continuous,recht2019tour, kumar2016optimal, levine2016end}, e.g., robotic in-hand manipulation~\cite{kumar2016optimal, levine2016end}.} 

However, there are only a few theoretical guarantees on PO methods as they often involve challenging non-convex optimization problems. To study their convergence and sample complexities, there has recently been a resurgent interest in PO methods for classical control problems \cite{fazel2018global,gravell2020learn,zhao2022sample,perdomo2021stabilizing, bu2019lqr, malik2019derivative, mohammadi2020linear, zhang2021policy,li2019distributed,zheng2021sample,bin2022towards}. For example, the seminal work \cite{fazel2018global} studies the well-known linear quadratic regulator (LQR) problem via PO methods. Though an optimal policy can be simply parameterized by a gain matrix, the quadratic cost is non-convex in the gain matrix space. A major contribution of \cite{fazel2018global} shows that the cost function is globally gradient dominated (aka Polyak-Lojasiewicz condition~\cite{polyak1963gradient, karimi2016linear}) with respect to (w.r.t.) the policy gain, which is indispensable to prove the global convergence of their PO methods.

%{This paper focuses on risk-constrained LQR problems~\cite{tsiamis2020risk,zhao2021infinitehorizon} which impose risk constraints on the system behavior. Such problems are originally proposed by \cite{tsiamis2020risk}, which introduces a new risk measure as a constraint, i.e., a one-step predicted state variability constraint. This formulation not only enables an explicit trade-off between control performance and the risk, but also yields a simple closed-form optimal controller. In our recent work~\cite{zhao2021infinitehorizon}, we extend the result to the infinite-horizon setting. In this paper, we solve the infinite-horizon risk-constrained problem~\cite{zhao2021infinitehorizon} in both model-based and sample-based settings, via PO methods rather than the dynamic programming (DP) methods used in  \cite{moore1997risk,ito2018risk,speyer1992optimal,borkar2014risk, chapman2019risk,tsiamis2020risk,zhao2021infinitehorizon}. }

%	
%	
%The Linear Quadratic Regulator (LQR) provides a fundamental optimal control framework for stochastic dynamical systems \cite{aastrom2012introduction,bertsekas1995dynamic}. 

Since the LQR problem only focuses on the quadratic regulation performance, the closed-loop system may be largely jeopardized by low-probability yet significant events, which is not allowed for safety-critical applications. To remedy it, risk-aware controllers have become natural choices \cite{moore1997risk,ito2018risk,speyer1992optimal, pan1996model,borkar2014risk, chapman2019risk,tsiamis2020risk,zhao2021infinitehorizon}. 
{In \cite{tsiamis2020risk}, a finite-horizon LQR problem with a variance-like constraint was first proposed, which is then extended to the infinite-horizon version in our previous work \cite{zhao2021infinitehorizon}. While both are solved via the model-based dynamic programming (DP), this paper studies the risk-constrained LQR (RC-LQR) problem of \cite{zhao2021infinitehorizon} under the PO framework in both model-based and sample-based settings.} In fact, various constrained LQ problems have also been studied via PO methods, e.g., the LEQG~\cite{zhang2021policy} and distributed LQG~\cite{li2019distributed}. 
%
%Under a quadratically invariant condition, the structural LQG problem~\cite{furieri2020learning} can be converted to an unconstrained one. In the PO for the LEQG~\cite{zhang2020policy},  the policy gradient methods have been addressed via an implicit regularization for the $\mathcal{H}_{\infty}$ constraint. Thus, the above LQR problems cannot be adapted to our setting. Overall, a primal-dual optimization landscape for our constrained LQRs is yet to be explored, which is our objective of this work.

%
% We solve it in both model-based and sample-based settings, via PO methods rather than the dynamic programming (DP) methods used in  \cite{moore1997risk,ito2018risk,speyer1992optimal,borkar2014risk, chapman2019risk,tsiamis2020risk,zhao2021infinitehorizon}. 
In sharp contrast to those PO works \cite{fazel2018global,gravell2020learn,zhao2022sample,perdomo2021stabilizing,bu2019lqr, malik2019derivative, mohammadi2020linear}, the non-convex variance-like constraint results in a fundamentally different optimization landscape. In particular, we lack the {\em global} gradient dominance property. 
Thus, a natural question is whether there still exists a good PO method that yields a globally optimal policy for the infinite-horizon RC-LQR problem. We provide a positive answer in this paper. {As the finite-horizon version \cite{tsiamis2020risk}, an optimal policy has also been shown in \cite{zhao2021infinitehorizon} to have an affine structure in the form of $u^*(x) = -K^*x + l^*$ with a gain matrix $K^*$ and a vector $l^*$}. We take this as a starting point, and propose here a novel primal-dual method where the primal and dual iterations alternatively compute an optimal policy-multiplier pair. 

%Though simple, it is yet challenging to establish theoretical guarantees for the primal-dual method, 
Even though the primal-dual method is conceptually simple, it is challenging to establish theoretical guarantees since (a) the optimization landscape of the Lagrangian function is yet unclear {(in fact, we only obtain that the Lagrangian under a fixed multiplier is {\em locally} gradient dominated, meaning that there may exist {\em multiple} optimal policies for the Lagrangian);} (b) the strong duality does not trivially hold in a non-convex constrained optimization problem (note that the strong duality is the key to primal-dual methods and is usually established for convex problems  \cite{bertsekas1997nonlinear}); and (c) exact gradients of both the Lagrangian and the dual function are unavailable. Our main contribution here lies in satisfactorily addressing the above issues, and further showing that the Lagrangian is also {\em coercive} with locally Lipschitz gradient, which along with local gradient dominance establishes the global convergence of our primal-dual method.

Clearly, the RC-LQR can be regarded as a special case of the long-studied constrained Markov decision problems (CMDPs)~\cite{altman1999constrained}. Strong duality for CMDPs has been proved, but only if the state-action space is finite~\cite{altman1999constrained} or the cost is uniformly bounded~\cite{ paternain2019constrained}, neither of which holds in the RC-LQR of this paper. To the best of our knowledge, we are the first to formally prove the strong duality for such a class of continuous CMDPs with quadratic costs.

Even though a similar policy gradient primal-dual framework has been adopted to solve  continuous CMDPs in \cite{ paternain2019safe, paternain2019constrained, bhatnagar2010actor, chow2017risk, tessler2018reward}, none of them can achieve global convergence. For example,  the primal-dual methods in  \cite{paternain2019safe, paternain2019constrained} have only been shown to converge to a neighborhood of the global optimum and can even lead to constraint violations. Even though it has been resolved in \cite{chow2017risk, tessler2018reward, bhatnagar2010actor},  their optimization landscape lends them resort to function approximations for optimal policies and thus can only achieve local convergence. In comparison, an optimal policy of our RC-LQR problem has an exact affine structure in the state feedback. While for finite CMDPs, the primal-dual methods are relatively easy and can ensure the convergence to a globally optimal policy~\cite{borkar2005actor, ding2020natural, ding2022convergence}.  It is worth mentioning that there are also other PO-based works~\cite{achiam2017constrained, pan2019risk,liu2020ipo,yu2019convergent} that do not follow a primal-dual framework, e.g.,  they leverage the interior-point method~\cite{liu2020ipo} and trust region method~\cite{yu2019convergent} to directly solve the constrained problem. Again, they still lack provable global convergence. 

The remainder of this paper is organized as follows. In Section \ref{sec:prob}, we formulate the infinite-horizon RC-LQR problem. In Section \ref{sec:primal-dual}, we approach it by proposing policy gradient primal-dual methods and recognizing the local gradient dominance property, based on which we prove the strong duality. In Section \ref{sec:mb} and Section \ref{sec:mf}, we propose primal-dual methods with convergence guarantees in model-based and sample-based settings, respectively. In Section \ref{sec:experiment}, we conduct simulations to validate our theoretical results. Concluding remarks of Section VII and five appendices complete the paper.
\section{Problem Formulation}\label{sec:prob}
Consider a discrete-time linear time-invariant stochastic system
\begin{equation}\label{equ:sys}
x_{t+1} = Ax_t + Bu_t+w_t,
\end{equation}
where $x_t\in\mathbb{R}^{n}$ and $u_t\in\mathbb{R}^{m}$ are the state and control vectors, and $\{w_t\}$ is an independently and identically distributed noise sequence.  

The infinite-horizon LQR problem aims to find a sequence of control policies $\{\pi_t\}$ to minimize a time-average cost, i.e.,
\begin{equation}\label{prob:lqr}
\begin{aligned}
\text {minimize}_{\{\pi_t\}} &~~ \lim\limits_{T \rightarrow  \infty} \frac{1}{T} \mathbb{E}\left[\sum_{t=0}^{T-1}(x_{t}^{\top} Q x_{t}+u_{t}^{\top} R u_{t})\right]\\
\text {subject to} &~~(\ref{equ:sys})~\text{and}~u_t = \pi_t(h_t, x_t),
\end{aligned}
\end{equation}
where $h_t = \{x_0,u_0,\cdots, x_{t-1},u_{t-1}\}$ is  the system history trajectory. The expectation is taken over the statistics of the noise sequence $\{w_t\}$. Throughout this paper, we make the following standard assumption~\cite{bertsekas1995dynamic}.
\begin{assum}
	\label{assumption}
	$Q$ is positive semi-definite and $R$ is positive definite. The pair $(A,B)$ is controllable and $(A,Q^{{1}/{2}})$ is observable.
\end{assum}

Under Assumption \ref{assumption}, solving (\ref{prob:lqr}) yields a unique optimal policy $\pi_t^*(x_t,h_t)=-Kx_t$ if the sequence $\{w_t\}$  has zero mean. Clearly, the LQR is risk-neutral as it only minimizes the quadratic cost, and the state may be substantially influenced by extreme noises, especially if $w_t$ has a heavy-tailed distribution. {To address it, the finite-horizon RC-LQR problem has been proposed in \cite{tsiamis2020risk}, which has been extended to the infinite-horizon case in our recent work \cite{zhao2021infinitehorizon}}, which has the following form 
\begin{equation}\label{equ:rclqr}
\begin{aligned}
&\text {minimize}_{\{\pi_t\}} ~ \lim\limits_{T \rightarrow  \infty} \frac{1}{T} \mathbb{E}\left[\sum_{t=0}^{T-1}(x_{t}^{\top} Q x_{t}+u_{t}^{\top} R u_{t})\right]\\
&\text {subject to} ~~(\ref{equ:sys}), ~u_t = \pi_t(h_t, x_t) ~\text{and} \\
&\lim\limits_{T \rightarrow  \infty} \frac{1}{T}  \mathbb{E}\left[ \sum_{t=0}^{T-1}(x_{t}^{\top} Q x_{t}-\mathbb{E}[x_{t}^{\top} Q x_{t}|h_t])^{2}\right] \leq \rho 
\end{aligned}
\end{equation}
where $\rho > 0$ is a user-defined constant to reflect our risk tolerance and $w_t$ has a finite 4th-order moment for tractability.  

{Different from \cite{tsiamis2020risk} and \cite{zhao2021infinitehorizon}, however, in this paper we re-solve the RC-LQR problem \eqref{equ:rclqr} via policy optimization (PO) methods}. 
{Following the notations of \cite{tsiamis2020risk},} let the mean and covariance of stationary noise $\{w_t\}$ be given by $\bar{w} = \mathbb{E}[w_t]$, $W = \mathbb{E}[(w_t - \bar{w})(w_t - \bar{w})^{\top}]>0$. Define
\begin{align}
M_{3} &= \mathbb{E}[(w_{t}-\bar{w})(w_{t}-\bar{w})^{\top} Q(w_{t}-\bar{w})], \\
m_{4} &= \mathbb{E}[(w_{t}-\bar{w})^{\top} Q(w_{t}-\bar{w})-\operatorname{tr}\{WQ\}]^{2}. \label{fouth}
\end{align}  
By  \cite[Theorem 1]{zhao2021infinitehorizon}, an optimal policy of (\ref{equ:rclqr}) has a time-invariant affine structure, i.e., $\pi_t^*(x_t,h_t)= -K^*x_t +l^*$, which also stabilizes the system \eqref{equ:sys} in the mean square sense. {Thus, there is no loss of optimality to solve \eqref{equ:rclqr} by focusing on the parameterized policies in the form of $u(x)=-Kx+l$, leading to the following optimization problem}
\begin{equation}\label{prob:new_rclqr}
\begin{aligned}
&\text{minimize}_{K,l} ~J(K,l) := \lim\limits_{T \rightarrow  \infty} \frac{1}{T} \mathbb{E} \left[ \sum_{t=0}^{T-1}(x_{t}^{\top} Q x_{t}+u_t^{\top} R u_t)\right] \\
&\text {subject to} ~~~(\ref{equ:sys}), u_t = -Kx_t + l, ~\text{and}~\\
&J_c(K,l):= \lim\limits_{T \rightarrow  \infty} \frac{1}{T}  \mathbb{E} \left[\sum_{t=0}^{T-1}(4x_{t}^{\top} QWQ x_{t} + 4x_{t}^{\top}QM_3  )\right] \leq \bar{\rho}
\end{aligned}
\end{equation}
with $\bar{\rho}= \rho - m_4 + 4\operatorname{tr}\{(WQ)^2\}$. 

The PO method for the risk-neutral LQR in (\ref{prob:lqr}) is shown to be globally convergent by random search~\cite{malik2019derivative, mohammadi2020linear} and policy gradient methods~\cite{fazel2018global,bu2019lqr}. However, the non-convex constraint in (\ref{prob:new_rclqr}) renders our problem much more involved and we resort to  the  duality theory to establish global convergence.  

\section{Primal-dual Methods for the Risk-constrained LQR}\label{sec:primal-dual}
In this section, we solve the RC-LQR problem (\ref{prob:new_rclqr}) via the primal-dual method.  {We first show that its Lagrangian function is coercive and locally gradient dominated.} Then, we establish strong duality. 
%
%Together with the coercivity and smoothness property, the global convergence of the primal-dual method is finally proved.
\subsection{Overview of Our Policy Gradient Primal-dual Method}
Let $X=[K, l]$ be the decision vector of (\ref{prob:new_rclqr}) and define the set of stabilizing policy by
$$
	\mathcal{S} = \{[K,l]|\rho(A-BK)<1, K \in \mathbb{R}^{m\times n}, l \in \mathbb{R}^m\}.
$$
Let $\mu \geq 0$ denote a Lagrange multiplier of (\ref{prob:new_rclqr}), $Q_{\mu} = Q+4 \mu Q W Q$ and $S = 2\mu QM_3$. Then, the Lagrangian is given as
\begin{equation}\label{def:L}
\begin{aligned}
\mathcal{L}(X,\mu) &= J(X)+\mu(J_c(X) - \bar{\rho}) \\
&= \lim\limits_{T \rightarrow  \infty} \frac{1}{T} \mathbb{E} \left[ \sum_{t=0}^{T-1}c_\mu(x_t, u_t)\right],
\end{aligned}
\end{equation}
where
\begin{equation}\label{def:lqtcost}
c_\mu(x_t, u_t)=x_{t}^{\top} Q_{\mu} x_{t}+ 2x_{t}^{\top}S +u_{t}^{\top} R u_{t}- \mu  \bar{\rho}
\end{equation}
is a reshaped cost with a non-negative weight $\mu$ to  balance the quadratic cost and the risk. Define the dual function as
\begin{equation}\label{def:D}
D(\mu) =  \min \limits_{X \in \mathcal{S}} \mathcal{L}(X,\mu).
\end{equation}

In the sequel, we refer to (\ref{prob:new_rclqr}) as the primal problem and 
\begin{equation}\label{prob:dual}
\mathop{\text{maximize}}\limits_{\mu \geq 0}~ D(\mu)
\end{equation}
as its dual problem.
Our primal-dual method is alternatively updated  as
\begin{subequations}
\label{primal-dual}
\begin{align}
X^{k}&\in \mathop{\text{argmin}} \limits_{X \in \mathcal{S}}~ \mathcal{L}(X,\mu^k), \label{prim_iterate} \\
\mu^{k+1}&= [\mu^k+ \zeta^k \cdot d^k]_+, \label{dual_iterate}
\end{align}
\end{subequations}
where the stepsize $\zeta^k>0$, $d^k$ is a subgradient of $D(\mu)$ at $\mu^k$ and $[x]_+=\max\{0, x\}$ for any $x\in\mathbb{R}$.

To achieve its global convergence, the strong duality property between the primal problem and the dual problem is essential.  Since (\ref{prob:new_rclqr}) is non-convex, it does not trivially hold. {Even though the Lagrangian in \eqref{def:L} is the LQR cost with a linear term,  its non-convex optimization landscape is yet unclear.} Thus, computing the primal update in \eqref{prim_iterate} is itself challenging.  In the rest of this section, we show that:  (a) $\mathcal{L}(X,\mu)$ is coercive over $\mathcal{S}$ and locally gradient dominated in Section \ref{subsec:property}, which is key to establish that a critical point of \eqref{prim_iterate} is globally optimal; (b) $\mathcal{L}(X,\mu)$ and its gradient are locally Lipschitz in Section \ref{subsec_local}, which implies a linear convergence rate of gradient methods for solving \eqref{prim_iterate}; (c) The strong duality property indeed holds in Section \ref{subsec_strong}. Combining these results prove the global convergence of \eqref{primal-dual}. Note that all the proofs on the properties of the Lagrangian are provided in Appendix \ref{Appendix:A} and \ref{apx:lip}. 

\subsection{Coercivity and Local Gradient Dominance of the Lagrangian}\label{subsec:property}

We first derive closed-form expressions for the Lagrangian and its gradient. For any $X \in \mathcal{S}$, the state of the system (\ref{equ:sys}) has a stationary distribution, the mean $\bar{x}_X$ and covariance $\Sigma_{K}$ of which satisfy 
\begin{align}\label{def:bar}
\bar{x}_X &= (A-BK)\bar{x}_X + Bl + \bar{w},\\
\label{def:sigma}
\Sigma_{K}&=W+(A-B K) \Sigma_{K}(A-B K)^{\top}.
\end{align}
%For technical proofs, we assume that $W > 0$, which can always be achieved by adding independent zero-mean Gaussian noises to $w_t$ without affecting the solution of (\ref{prob:new_rclqr}).
%It can be easily observed that $\mathcal{L}(X,\mu)$ is finite if and only if
Then, we define the value function under $X$ associated with the reshaped cost $c_\mu(x_t,u_t)$ as
$$
V_{X}(x)=\mathbb{E}\left[\sum_{t = 0}^{\infty} (c_\mu(x_{t}, u_{t}) -\mathcal{L}(X,\mu)) \bigg| x_{0}=x \right],
$$
where $\mathbb{E}[\cdot]$ takes expectation under a fixed policy $X\in\mathcal{S}$.
%which is finite~\cite{yang2019provably}.
Moreover, let $P_K \geq 0$ satisfy the following Lyapunov equation
\begin{equation}\label{equ:riccati}
	P_{K} = Q_{\mu} +K^{\top} R K + (A-B K)^{\top} P_{K}(A-B K)
\end{equation}
and define
\begin{align*}
&E_{K}=(R+B^{\top} P_{K} B) K-B^{\top} P_{K} A, \\
&R_K = R + B^{\top} P_K B, ~~
V = (I - (A-BK))^{-1}.
\end{align*}

We show that $V_{X}(x)$ is quadratic and provide a closed-form of $\mathcal{L}(X,\mu)$.
\begin{lemma}\label{prop:lag}
	For any $X \in \mathcal{S}$, it follows that
	\begin{align*}\label{equ:value}
	\text{(\romannumeral1)}~~~~ V_{X}(x)&= x^{\top}P_Kx + g_{X}^{\top}x + z_{X}, \\
	\text{(\romannumeral2)}~~ \mathcal{L}(X,\mu)& = \mathrm{tr}\{P_K(W + (Bl+\bar{w})(Bl+\bar{w})^{\top})\}\\
	&+ g_{X}^{\top}(Bl+\bar{w})  + l^{\top}Rl - \mu \bar{\rho},
	\end{align*}
	where $g_{X}^{\top} = 2(-l^{\top}E_K+S^{\top} + \bar{w}^{\top}P_K(A-BK))V$
	and $z_{X}$ is a constant irrespective of $x$.
\end{lemma}

Moreover, the gradient of $\mathcal{L}(X,\mu)$ w.r.t. $X$ is explicitly given in the following lemma.
\begin{lemma}\label{prop:grad}
	For any $X \in \mathcal{S}$, the gradient of $\mathcal{L}(X,\mu)$ in $X$ is
	\begin{equation}\label{equ:gradient}
	\nabla_X \mathcal{L}(X,\mu) = 2[
	E_K,G_X
	] \Phi_X,
	\end{equation}
	where $G_X = R_Kl + B^{\top}P_K\bar{w} + \frac{1}{2} B^{\top}g_X$ and $\Phi_X$ is an ergodic matrix
	\begin{equation}\label{def:phi}
	\begin{aligned}
	\Phi_X &= \lim\limits_{T \rightarrow  \infty} \frac{1}{T}\mathbb{E}\left[\sum_{t=0}^{T-1} 		
	\begin{bmatrix}
	x_t \\
	-1
	\end{bmatrix}
	\begin{bmatrix}
	x_t \\
	-1
	\end{bmatrix}^{\top}\right]\\
	&=
	\begin{bmatrix}
	\Sigma_K + \bar{x}_{X}\bar{x}_{X}^{\top} &  -\bar{x}_{X}\\
	-\bar{x}_{X}^{\top} & 1
	\end{bmatrix}>0.
	\end{aligned}	
	\end{equation}
\end{lemma}

Since $\Phi_X>0$,  letting $\nabla_X \mathcal{L}(X,\mu)=0$ yields a unique critical point 
\begin{equation}\label{equ:stationary}
X^*(\mu)= [K_\mu,l_\mu]
\end{equation}
with $K_\mu = R_{K_\mu}^{-1}B^{\top} P_{K_\mu} A$ and $l_\mu = -R_{K_\mu}^{-1}B^{\top}V^{\top}(P_{K_\mu}\bar{w}+S)$. 
%We study the properties of $\mathcal{L}(K,l,\mu)$ to facilitate the optimization in (\ref{def:D}).
%
%The minimization of $\mathcal{L}(X,\mu)$ in (\ref{def:D}) is a non-convex optimization problem as both the objective function and the stabilizing policy set $\mathcal{S}$ are non-convex, which poses challenges in solving (\ref{def:D}) with standard gradient-based methods. In the PO for classical LQR problems (\ref{prob:lqr})~\cite{fazel2018global}, this is alleviated by observing that the objective function is globally gradient dominated.
%we refer to as the {\em global} gradient dominance and
Now,  we are ready to show two important properties of the Lagrangian.

\begin{lemma}[{\bf Coercivity}]\label{lem:coer}
	 Under a fixed $\mu>0$, 
	$\mathcal{L}(X,\mu)$ is coercive in $X$ in the sense that
	$
	\lim_{X \rightarrow \partial \mathcal{S}} \mathcal{L}(X,\mu) = +\infty,
	$
	where $\partial \mathcal{S}$ denotes the boundary of $\mathcal{S}$, and has a compact $\alpha$-sublevel set 
	\begin{equation}\label{def:sublevel}
	\mathcal{S}_{\alpha}  = \{ X |\mathcal{L}(X,\mu) \leq \alpha \}.
	\end{equation}

\end{lemma}

 \begin{defi}For a differentiable function $f(x): \mathbb{R}^n \rightarrow \mathbb{R}$ with a finite global minimum $f^*$, it is  gradient dominated over a set $\mathcal{X}\subseteq \text{dom}(f) $ if
\begin{equation}\label{dominated}
f(x) - f^* \leq \lambda_{\mathcal{X}} \| \nabla f(x) \|^2, ~~\forall x \in \mathcal{X}~\text{for some}~ \lambda_\mathcal{X} > 0.
\end{equation}
\end{defi} 

If $\mathcal{X}=\text{dom}(f)$, it reduces to the Polyak-Lojasiewicz condition ~\cite{polyak1963gradient, karimi2016linear} which is key to the global convergence of \cite{fazel2018global}.   In this paper, we can only show that \eqref{dominated} holds for some proper subset of $\text{dom}(f)$, and for distinction refer to them as the {\em global} and {\em local}  gradient dominance, respectively. 

For a given policy $X\in\mathcal{S}$, define a truncated value function
$$
V_{X}^{T}(x)= \mathbb{E}\left[\sum_{t = 0}^{T-1} \left(c_\mu(x_{t}, u_{t}) -\mathcal{L}(X,\mu)\right) \bigg| x_{0}=x\right]
$$
and an advantage function {
$
A_{X}^{T}(x, u)= c_{\mu}(x,u)-\mathcal{L}(X,\mu) + \mathbb{E}[V_{X}^{T}(x_{t+1})|x_t=x, u_t=u]-V_{X}^{T}(x)
$. 
}
Then, the Lagrangian difference between the two stabilizing policies can be described by the advantage function.
\begin{lemma}\label{lem:difference}
	Let $\{x_t'\}$ and $\{u_t'\}$ be sequences generated by the stabilizing policy $X'\in \mathcal{S}$. For any $X\in\mathcal{S}$, it follows that
	\begin{align*}
	&\text{(\romannumeral1)}~\mathcal{L}(X',\mu)-\mathcal{L}(X,\mu)=-
	\lim\limits_{T \rightarrow  \infty} \frac{1}{T}\mathbb{E}\left[\sum_{t=0}^{T-1}A_{X}^{T}(x_t', u_t')\right]\\
	&\text{(\romannumeral2)}~\lim\limits_{T \rightarrow  \infty} A_{X}^{T}(x, -K'x+l') \\
	&~~~~= ((K'-K)x-(l'-l) + R_K^{-1}(E_Kx-G_X))^{\top}\\
	&~~~~~~~\times R_K ((K'-K)x-(l'-l) + R_K^{-1}(E_Kx-G_X))\\
	&~~~~~~~- (E_Kx-G_X)^{\top}R_K^{-1}(E_Kx-G_X), ~\forall x \in \mathbb{R}^n.
	\end{align*}		
\end{lemma}

Lemma \ref{lem:difference} is consistent with Lemma 10 in \cite{fazel2018global}, though we focus  on the ergodic cost here. 
\begin{lemma} \label{lem:gradient dominance}
	For any $X \in \mathcal{S}$, it  holds that
	\begin{equation}\label{equ:lemma2}
	\mathcal{L}(X,\mu)-\mathcal{L}^*(\mu)\leq \frac{\|\Phi^*\|}{4\underline\sigma(R)\underline\sigma(\Phi_X)^2}\cdot\operatorname{tr} \{
	\nabla_X\mathcal{L}^{\top}\nabla_X\mathcal{L}	
	\},
	\end{equation}
	where $\nabla_X \mathcal{L} = \nabla_X \mathcal{L}(X,\mu)$ is given in (\ref{equ:gradient}), $\Phi^*$ is the ergodic matrix (\ref{def:phi}) under $X^*(\mu)$ of (\ref{equ:stationary}), $\underline{\sigma}(\cdot)$ returns the minimum eigenvalue of a positive definite matrix, and $\mathcal{L}^*(\mu) = \min_{X\in S}\mathcal{L}(X,\mu)$.
\end{lemma}

{Since $\lim_{l \rightarrow \infty}\underline\sigma(\Phi_X)=0$ (cf. \eqref{def:bar} and \eqref{def:phi}), the coefficient on the right hand side of \eqref{equ:lemma2} is unbounded, in contrast to the case of the LQR~\cite{fazel2018global}, where it is a finite constant, i.e., their quadratic cost is globally gradient dominated.} The good news here is that it is also finite over the $\alpha$-sublevel set $\mathcal{S}_{\alpha}$ in \eqref{def:sublevel}, as established below. 
\begin{lemma}[{\bf Local gradient dominance}] \label{lem:local_grad} For any $\mu>0$,
	$\mathcal{L}(X, \mu)$ is gradient dominated over its $\alpha$-sublevel set, i.e., 
	$$
	\mathcal{L}(X,\mu)-\mathcal{L}^*(\mu)
	\leq \lambda_{\alpha}\cdot\operatorname{tr} \{
	\nabla_X\mathcal{L}^{\top}\nabla_X\mathcal{L}	
	\}, \forall X\in\mathcal{S}_{\alpha},
	$$
	where $\lambda_{\alpha} = {\|\Phi^*\|}/({4\underline\sigma(R) \cdot \sigma_{\alpha}^2}) >0$ is a constant over $\mathcal{S}_{\alpha}$ in \eqref{def:sublevel} and $\sigma_{\alpha} = \min_{X \in \mathcal{S}_{\alpha}}  \underline\sigma(\Phi_X)>0$.
\end{lemma}
\begin{proof}
	Since $\Phi_X > 0$ is continuous in $X$,  then $\underline\sigma(\Phi_X)$ can be lower bounded by a positive constant over the compact set $\mathcal{S}_{\alpha}$. The result then follows.
\end{proof}

{Since Lemma \ref{lem:local_grad} holds for any $\alpha>0$, joint use of  coercivity is sufficient for finding a global minimizer of \eqref{prim_iterate}.}
\begin{theo}\label{theorem:unique} For any $\mu>0$, the critical point $X^*(\mu)$ in (\ref{equ:stationary}) is the unique global minimizer of $\mathcal{L}(X,\mu)$.
\end{theo}
\begin{proof} It is straightforward from \eqref{equ:stationary} and Lemma \ref{lem:local_grad}.
\end{proof}
%That is, the primal iteration in \eqref{prim_iterate} can be trivially solved, i.e., $X^k = X^*(\mu^k)$.
%In fact, this local property would be in effect \textit{globally} if the optimization sequence $\{X_i\}$ is guaranteed to be contained in $\mathcal{S}_{\alpha}$. We leverage this observation to develop our sample-based primal-dual algorithms.
\subsection{Locally Lipschitz  Gradient of the Lagrangian}\label{subsec_local}
%If  $f(x)$ is both globally gradient dominated and smooth, then a gradient-based method converges at a linear rate~\cite{malik2019derivative, karimi2016linear}. Unfortunately, the  smoothness of the Lagrangian does not hold here. For the standard LQR, this problem is addressed by the so-called \textit{almost smoothness} property~\cite{fazel2018global}.  

For a fixed $\mu$,  we show in this subsection that both $\mathcal{L}(X, \mu)$ and its gradient  are locally Lipschitz continuous.

%Moreover, its gradient also has a local Lipschitz constant which can be utilized to bound the stepsize of gradient methods.

\begin{lemma} \label{lem:cost_diff}
	 For any pair of stabilizing policies $X$ and $X'$, the gap of their Lagrangians  is given as
	\begin{align*}
	&\mathcal{L}(X',\mu)-\mathcal{L}(X,\mu)\\
	&= \operatorname{tr} \{(K'-K)^{\top}R_K(K'-K) (\Sigma_{K'} + \bar{x}'\bar{x}'^{\top}) \\
	&+2 (K'-K)^{\top} E_{K} (\Sigma_{K'} + \bar{x}'\bar{x}'^{\top})-2G_X^{\top}(K'-K)\bar{x}'\\
	&- 2(l'-l)^{\top}R_K(K'-K)\bar{x}'+(l'-l)^{\top}R_K(l'-l)\\
	&-2(l'-l)^{\top}E_K\bar{x}'+2(l'-l)^{\top}G_X\}
	\end{align*}
	where $\bar{x}'$ denotes $\bar{x}_{X'}$ in \eqref{def:bar} for notational simplicity.
\end{lemma}

\begin{lemma}[{\bf Locally Lipschitz Lagrangian and gradient}] \label{lem:lip} For any $X\in\mathcal{S}$, there exist positive scalars $(\xi_{X}, \beta_X, \gamma_X)$ such that for any  $X'\in\mathcal{S}$ and $\| X'-X\| \leq \gamma_X$, it holds 
\begin{align*}
	&|\mathcal{L}(X',\mu) - \mathcal{L}(X,\mu)|   \leq \xi_X \| X'-X\|, ~\text{and} \\
	&\|\nabla_X \mathcal{L}(X',\mu) - \nabla_X \mathcal{L}(X,\mu)\|   \leq \beta_X \| X'-X\|.
	\end{align*}
\end{lemma}

The scalars $(\xi_{X}, \beta_X, \gamma_X)$ are polynomials of $\|A\|$, $\|B\|$, $\underline{\sigma}(Q_{\mu})$, $\underline{\sigma}(R)$, and are uniformly bounded over a compact set. 

{Comparing with \cite{fazel2018global}, i.e., letting $\mu$ be fixed and $M_3=0,\bar{w}=0$, we have $S=0$ and $l=0$ in (\ref{def:lqtcost}) and the Lagrangian reduces to the standard LQR cost. Then, Section \ref{subsec:property} and \ref{subsec_local} recover the results in \cite[Section 3]{fazel2018global}.}

\subsection{Strong Duality}\label{subsec_strong}
{In this subsection, we show that the strong duality between the primal problem  (\ref{prob:new_rclqr}) and dual problem \eqref{prob:dual} holds.}
\begin{lemma}\label{lem:continuity}
	 Both the policy $X^*(\mu)$ in \eqref{equ:stationary} and the constraint function $J_c(X^*(\mu))$ are continuous over $\mu\in[0,\infty)$.
\end{lemma}
\begin{proof}
For $\mu \geq 0$ and $X \in \mathcal{S}$, it follows from  (\ref{equ:riccati}) that the Lyapunov equation yields a unique $P_K >0$, which jointly with (\ref{equ:stationary}) implies that $X^*(\mu)$ is continuous in $\mu \ge 0$.  The continuity of $J_c(X^*(\mu))$ can be established by using the arguments in \cite[Lemma 3.6]{bu2019lqr}.
\end{proof}

	{Note that the continuity in Lemma \ref{lem:continuity} is a strong result and usually lacks in the primal-dual framework.  Particularly, it holds only if $X^*(\mu)$ in \eqref{primal-dual}  is unique, which is not the case for a  general non-convex optimization problem. We now formally prove the strong duality result under Slater's condition, which essentially follows from \cite[Theorem 3]{tsiamis2020risk} and \cite[Chapter 6]{bertsekas1997nonlinear}.}
	 \begin{assum}[{\bf Slater's condition}]\label{assum:slater}
	There exists a policy $\widetilde{X} \in \mathcal{S}$ such that $J_c(\widetilde{X})<\bar{\rho}$. 
	 \end{assum}
\begin{theo}[{\bf Strong duality}]\label{theorem:duality}
	Under Assumption \ref{assum:slater}, there is no duality gap between the primal problem (\ref{prob:new_rclqr}) and the dual problem \eqref{prob:dual}.
\end{theo}
\begin{proof}
%By \cite[Proposition 5.1.5]{bertsekas1997nonlinear}, a pair of solutions $(X_0,\mu_0)$ are optimal solutions to the primal and dual problems (\ref{prob:new_rclqr}) and (\ref{prob:dual}) with zero duality gap {\em if and only if} the following optimality conditions hold, i.e.,
%\begin{equation}\label{equ:saddle}
%\begin{aligned}
%\mathcal{L}(X^*,\mu^*) & = \min_{X \in \mathcal{S}}\mathcal{L}(X,\mu^*), \\
%J_c(X^*) & \leq \bar{\rho}, \\
%\mu^*(J_c(X^*)- \bar{\rho})&=0.
%\end{aligned}
%\end{equation}
Define
\begin{equation}\label{def:mu}
\mu^{*} = \inf \{\mu \geq 0| J_c(X^*(\mu)) \leq \bar{\rho}\},
\end{equation}
where $X^*(\mu) \in \mathop{\text{argmin}}_{X \in \mathcal{S}} \mathcal{L}(X,\mu)$. 
By \cite[Proposition 6.1.5]{bertsekas1997nonlinear}, it is sufficient to show that (a) $\mu^*$ is finite, and (b) the policy-multiplier pair $(X^*, \mu^*)$ with $X^*=X^*(\mu^*)$  satisfies the following optimality conditions 
\begin{equation}\label{equ:saddle}
\begin{aligned}
\mathcal{L}(X^*,\mu^*)  &= \min_{X \in \mathcal{S}}\mathcal{L}(X,\mu^*), \\
J_c(X^*)  &\leq \bar{\rho}, \\
\mu^*(J_c(X^*)- \bar{\rho})&=0.
\end{aligned}
\end{equation}
(a) By Assumption \ref{assum:slater}, there exists a constant $a>0$ such that $J_c(\widetilde{X})+a \leq \bar{\rho}$. We prove by contradiction and assume that for all $\mu \geq 0$, $J_c(X^*(\mu)) > \bar{\rho}$. Then, 
\begin{align*}
J(\widetilde{X}) &\geq D(\mu) - \mu (J_c(\widetilde{X}) - \bar{\rho})\\
&= J(X^*(\mu)) + \mu (J_c(X^*(\mu)) - J_c(\widetilde{X})) \\
&\geq J(X^*(\mu)) + \mu (J_c(X^*(\mu)) - \bar{\rho} + a) \\
&> J(X^*(\mu)) + \mu a.
\end{align*}
Letting $\mu \rightarrow \infty$ implies that $J(\widetilde{X}) > \infty$, which contradicts Slater's condition that $\widetilde{X} \in \mathcal{S}$. Thus, $\mu^*$ in (\ref{def:mu}) is finite.

(b) Clearly, we only need to verify that $\mu^*(J_c(X^*)- \bar{\rho})=0$. If $\mu^* = 0$, then it trivially holds. If $\mu^* >0$, it follows from \eqref{def:mu} that $J_c(X^*(0)) > \bar{\rho}$. Since $\mu^*$ is finite, there must exist a $\mu'>0$ such that $J_c(X^*(\mu')) \leq \bar{\rho}$. The continuity of  $J_c(X^*(\mu))$ in Lemma \ref{lem:continuity} implies that $J_c(X^*)=\bar{\rho}$.
\end{proof}

%	{Though this paper considers a single constraint $J_c(X) \leq \bar{\rho}$, the strong duality result in Theorem \ref{theorem:duality} can be extended for the general case with multiple constraints in quadratic form, by using mathematical induction. }

%The policy gradient primal-dual method for solving (\ref{prob:new_rclqr}) can be decomposed into two alternative phases. Firstly, we minimize $\mathcal{L}(X,\mu)$ under the current multiplier $\mu$ via policy gradient methods. Then, we compute the subgradient of $D(\mu)$ and update the multiplier towards the direction of maximizing $D(\mu)$. In the sequel, we develop two primal-dual algorithms to solve (\ref{prob:new_rclqr}) in both model-based and sample-based settings.

\section{Policy Gradient Primal-dual Algorithm for the Model-based Setting}\label{sec:mb}
In the model-based setting, we assume that all the parameters in (\ref{equ:sys}) is known and propose three gradient-based methods with linear convergence to solve (\ref{prim_iterate}). Then, we develop a primal-dual method in the form of  (\ref{primal-dual})  with global convergence to solve (\ref{prob:new_rclqr}).
\subsection{Policy Gradient Methods for Solving (\ref{prim_iterate})}\label{subsec}
To solve (\ref{prim_iterate}), we consider three widely-used policy gradient methods~\cite{fazel2018global,bu2019lqr,gravell2020learn}. Let $X'$ be the one-step updated policy and $\eta$ be the stepsize. The update rules are given by
\begin{equation}\label{def:pgmethod}
	\hspace{-0.07cm}\begin{aligned}
	\textbf{Policy Gradient (PG):} ~X' &=X-\eta \nabla_X \mathcal{L}\\
	\textbf{Natural PG (NPG):} ~X'&=X-\eta \nabla_X \mathcal{L} \cdot \Phi_X^{-1} \\
	\textbf{Gauss-Newton (GN):} ~X' &=X-\eta R_K^{-1}\cdot \nabla_X \mathcal{L} \cdot \Phi_X^{-1}\\
	\end{aligned}
\end{equation}
where $\nabla_X \mathcal{L}$ and $\Phi_X$ can be computed via \eqref{equ:gradient} and \eqref{def:phi}, respectively. The NPG update is related to the gradient over a Riemannian manifold, while the GN update is one type of quasi-Newton update. 

For simplicity, we follow  \cite{fazel2018global,gravell2020learn, bu2019lqr, malik2019derivative, mohammadi2020linear, zhang2021policy,li2019distributed,zheng2021sample} to assume the access of an initial stabilizing policy $X^{(0)}\in\mathcal{S}$. Note that this can be relaxed via the PO methods; see e.g., \cite{perdomo2021stabilizing,zhao2022sample}.  

The key to the linear convergence of \eqref{def:pgmethod} is to find an appropriate stepsize such that  (\ref{def:pgmethod}) yields a stabilizing $X'$ and decreases the Lagrangian per iteration, which is formally stated below. Note that the  proof is given in Appendix \ref{apx:pg}.
%Denote by $\Phi'$ the ergodic matrix of $X'$ in \eqref{def:phi}. %In particular, the GN method enjoys a super linear rate when $\eta = \frac{1}{2}$.
\begin{theo}\label{theo:pg}Define the compact sublevel set $$\mathcal{S}_0 = \{ X |\mathcal{L}(X,\mu) \leq \mathcal{L}(X^{(0)},\mu)\}$$ and $\sigma_0 = \min_{X \in \mathcal{S}_0} \underline\sigma(\Phi_X)$.   
If $X\in\mathcal{S}_0$ and $\eta$ in \eqref{def:pgmethod} is appropriately selected, then there exists a finite $\beta\in(0,1)$ such that
$$
\mathcal{L}(X',\mu) - \mathcal{L}^*(\mu) \leq \left(1 - \beta\right)(\mathcal{L}(X,\mu) - \mathcal{L}^*(\mu)).
$$

Moreover, (a) $0<\eta\leq 1/2$ and $\beta= {2\eta \sigma_0}/{\|\Phi^*\|}$ for the GN update; (b) 
$0<\eta\leq {1}/({2\|R_{K_0}\|})$ and $\beta={ 2\eta \sigma_0 \underline{\sigma}(R)}/{\|\Phi^*\|}$ for the NPG update; and (c) $\eta$ is a polynomial in problem parameters and $\beta={ 2\eta \sigma_0^2 \underline{\sigma}(R)}/{\|\Phi^*\|}$ for the PG update. 
\end{theo}
%\begin{theo}\label{theorem:gn}
%	If $0<\eta\leq 1/2$ and $X\in\mathcal{S}_0$, the GN update
%	\begin{equation}\label{equ:GNupdate}
%		X' =X-2 \eta R_K^{-1} [E_K~~G_X]
%	\end{equation}
%	converges to $X^*(\mu)$ at a linear rate, i.e.,
%	$$
%	\mathcal{L}(X',\mu) - \mathcal{L}^*(\mu) \leq \left(1 - {2\eta \sigma_0}/{\|\Phi^*\|}\right)(\mathcal{L}(X,\mu) - \mathcal{L}^*(\mu)).
%	$$
%\end{theo}
%
%
%\begin{theo}\label{theorem:ng}
%	If $0<\eta\leq {1}/({2\|R_{K_0}\|})$ and $X\in\mathcal{S}_0$, the NPG update
%	\begin{equation}\label{equ:ng_update}
%	X' =X-2 \eta [E_K~~G_X]
%	\end{equation}
%	converges to $X^*(\mu)$ at a linear rate, i.e.,
%	$$
%	\mathcal{L}(X',\mu) - \mathcal{L}^*(\mu) \leq \left(1 - { 2\eta \sigma_0 \underline{\sigma}(R)}/{\|\Phi^*\|}\right)(\mathcal{L}(X,\mu) - \mathcal{L}^*(\mu)).
%	$$
%\end{theo}
%\begin{theo}\label{theorem:pg}
%	For an appropriate stepsize $\eta$ that is polynomial in problem parameters, e.g., $\|A\|$, $\|B\|$, ${\underline{\sigma}(W)}$, ${\underline{\sigma}(Q_{\mu})}$, ${\underline{\sigma}(R)}$ and $X\in\mathcal{S}_0$, the gradient update
%	\begin{align*}
%	X' =X-2 \eta [E_K~~G_X]\Phi_X
%	\end{align*}
%	converges to  $X^*(\mu)$ at a linear rate, i.e.,
%	$$
%	\mathcal{L}(X',\mu) - \mathcal{L}^*(\mu) \leq \left(1 - { 2\eta \sigma_0^2 \underline{\sigma}(R)}/{\|\Phi^*\|}\right)(\mathcal{L}(X,\mu) - \mathcal{L}^*(\mu)).
%	$$
%\end{theo}

\begin{table}[t]
	\caption{{Comparison of three gradient methods (I-best, III-worst).}}
	\label{table}
	\begin{center}
		\begin{tabular}{|c|c|c|c|}
			\hline
			 & convergence rate &stepsize & complexity\\
			\hline
			PG & III & III  & III \\
					
			 NPG& II & II &I   \\
				
			 GN& I & I & II \\			
			\hline
		\end{tabular}
		\label{tab1}
	\end{center}
\end{table}

{We provide a comparison of the three methods of \eqref{def:pgmethod} in Table \ref{table}.} Since the NPG and GN updates use more information, e.g., $\Phi_X$ and $R_K$, they tend to use less conservative stepsizes and achieve better convergence rates. {
Even though the PG update is given in the simplest form in \eqref{def:pgmethod},  the updates of $K$ and $l$  in GN and NPG can be decoupled, e.g., the NPG is rewritten as 
\begin{align*}
K' = K - 2\eta E_K,~\text{and}~~ l' = l - 2 \eta G_X
\end{align*}
which can reduce the  computational complexity per update. Nonetheless, their computational complexities are essentially the same as $\mathcal{O}(n^3)$.} Interestingly, the GN update with stepsize $\eta = 1/2$ is equivalent to the policy iteration and achieves a superlinear convergence rate~\cite{bu2019lqr} which is also confirmed via simulation in Section \ref{sec:experiment}. 
%This implies that the computation burden at each step can be reduced by starting optimizing $l$ after $K$ has converged, with a sacrifice of overall runtime. We leave it to our future work.
 	
\subsection{A Model-based Primal-dual Algorithm}
\begin{algorithm}[t!]
	\caption{The model-based primal-dual algorithm for the risk-constrained LQR}
	\label{alg:model-based}
	\begin{algorithmic}[1]
		\Require A randomly initialized multiplier $\mu^1\geq 0$, and a set of stepsizes $\{\zeta^k\}$.
		\For{$k=1,2,\dots$}
		%\State \textbf{Step 1: Minimize the Lagrangian}
		\State Solve $X^{k} = \text{argmin}_{X \in \mathcal{S}}~ \mathcal{L}(X,\mu^k)$ via (\ref{def:pgmethod}).
		%\State \textbf{Step 2: Dual ascent}
		\State Compute a subgradient $d^k$ by (\ref{def:sub}) and Lemma \ref{lem:risk_expression}.
		\State Update the multiplier by $\mu^{k+1} = [\mu^{k} + \zeta^k \cdot d^k]_{+}$.
		\EndFor	
		%\Ensure Policy-multiplier pair $(X^M,\mu^M)$.
	\end{algorithmic}
\end{algorithm}
By duality theory~\cite{nesterov2013introductory,nedic2009subgradient}, a subgradient in (\ref{dual_iterate}) is 
\begin{equation}\label{def:sub}
d^k = J_c(X^k) - \bar{\rho},
\end{equation}
where $X^k$ is given in (\ref{prim_iterate}) and $J_c(X^k)$ is computed by the following lemma.
\begin{lemma}\label{lem:risk_expression}
	For a stabilizing policy $X \in \mathcal{S}$, we have
	$$
	J_c(X) = \mathrm{tr}\{P_c(W + (Bl+\bar{w})(Bl+\bar{w})^{\top})\} + g_c^{\top}(Bl+\bar{w}),
	$$
	where $P_c > 0$ is a unique solution of the Lyapunov equation
	$$
	P_{c} = 4QWQ + (A-B K)^{\top} P_{c}(A-B K),
	$$
	and
	$ g_{c}^{\top} = 2((Bl+\bar{w})^{\top}P_c(A-BK)+ 2M_3^{\top}Q)V.$
\end{lemma}
\begin{proof}
	The proof is similar to that of Lemma \ref{prop:lag}.
\end{proof}

Our model-based primal-dual method is summarized in Algorithm \ref{alg:model-based}. In general, the primal iteration will not converge to a feasible solution unless the subdifferential of the dual function is a singleton~\cite{bertsekas1997nonlinear, boyd2004convex}.  Fortunately,  Theorem \ref{theorem:unique} implies that $X^k$ is the unique minimizer of $\mathcal{L}(X, \mu^k)$. Since $X^k$ is always able to stabilize the system, the subgradient (actually gradient) $d^k$ and  $\mu^{k}$ are {\em uniformly bounded}. Jointly with the concavity of $D(\mu)$, it follows from \cite[Theorem 3]{zhao2021infinitehorizon} that  Algorithm \ref{alg:model-based}  converges globally.
\begin{theo}\label{theo:model-based}
Let $\bar{\mu}^k = \frac{1}{k} \sum_{i=1}^{k} \mu^{i}$ and $\zeta^k = \mathcal{O}(k^{-1/2})$. Under Assumption \ref{assum:slater}, Algorithm \ref{alg:model-based} yields
$$
  D^*-	D(\bar{\mu}^k) \leq \mathcal{O}(k^{-1/2}),
$$
where the maximum of the dual function $D^* = \max_{\mu \geq 0}~ D(\mu)$ is finite. 

\end{theo}
%\begin{proof}
%	By the definition of projection, it follows that
%	\begin{align*}
%		\|\mu^{i+1} - \mu^*\|^2 &\leq \| \mu^{i} - \mu^* + \zeta^i \cdot d^i\|^2 \\
%		&= \|\mu^{i} - \mu^*\|^2 +2\zeta^i d^{i{\top}}(\mu^{i} - \mu^*) + (\zeta^i)^2\|d^i\|^2\\
%		&\leq \|\mu^{i} - \mu^*\|^2 +2\zeta^i(D(\mu^{i}) - D^*) + (\zeta^i)^2b^2.
%	\end{align*}
%	
%	Then, rearranging it yields that
%	\begin{align*}
%		D^*-D(\mu^{i}) \leq \frac{\|\mu^{i} - \mu^*\|^2}{2\zeta^i} -\frac{\|\mu^{i+1} - \mu^*\|^2}{2\zeta^i} + \frac{\zeta^i b^2}{2}.
%	\end{align*}
%	
%	Summing up from $i=1$ to $k$ and noting $\zeta^i \geq \zeta^{i+1}$, it follows that
%	\begin{align*}
%		&\sum_{i=1}^{k}(D^*-D(\mu^{i})) \leq -\frac{1}{2\zeta^{k+1}}\|\mu^{k+1}-\mu^*\| + \frac{b^2}{2} \sum_{i=1}^{k} \zeta^i \\
%		&~~~+ \frac{1}{{2\zeta^1}}\|\mu^{1} - \mu^*\|^2 +  \frac{1}{2}\sum_{i=1}^{k}(\frac{1}{\zeta^{i+1}} - \frac{1}{\zeta^i})\|\mu^{i+1}-\mu^*\|^2 \\
%		&\leq \frac{2}{\zeta^k}e^2 + \frac{b^2}{2} \sum_{i=1}^{k} \zeta^i.
%	\end{align*}
%
%	By Jenson's inequality, one can easily obtain that
%	\begin{align*}
%		D^*-	D(\bar{\mu}^k) \leq \frac{2}{k\zeta^k}e^2 + \frac{b^2}{2k} \sum_{i=1}^{k} \zeta^i.
%	\end{align*}
%	
%	The proof follows by noting that $\zeta^i = \frac{1}{be} \sqrt \frac{2}{i}$.
%\end{proof}
\begin{proof}
	It is similar to that of \cite[Theorem 3]{zhao2021infinitehorizon} and omitted for saving space.
\end{proof}

{By Lemma \ref{lem:risk_expression}, a simple bisection method could be adopted to solve $\mu^*$ in \eqref{def:mu} for the model-based setting. However,  it is unclear how to adopt it for the sample-based setting as the constraint function can only be randomly evaluated as well.}

\section{Policy Gradient Primal-dual Algorithm for the Sample-based Setting} \label{sec:mf}
If $(A,B)$ in \eqref{equ:sys} is unknown, both $\nabla_X \mathcal{L}(X,\mu)$ in (\ref{equ:gradient}) and $d^k$ in \eqref{def:sub} cannot be computed directly. In the sample-based setting, we estimate them via  system trajectories and develop a sampled-based primal-dual algorithm with global convergence.

Specifically,  assume that there is an oracle to return noisy values of $\mathcal{L}(X,\mu)$ and $J_c(X)$ viz
\begin{equation}\label{oracle}
\begin{aligned}
&\widehat{\mathcal{L}}(X,\mu) =  \lim\limits_{T \rightarrow  \infty} \frac{1}{T} \sum_{t=0}^{T-1}c_\mu(x_t, u_t),  ~~ \text{and} \\
&\widehat{J}_c(X)  = \lim\limits_{T \rightarrow  \infty} \frac{1}{T}  \sum_{t=0}^{T-1}(4x_{t}^{\top} QWQ x_{t} + 4x_{t}^{\top}QM_3),
\end{aligned}
\end{equation}
 where $\{x_t\}$ and $\{u_t\}$ denote the states and control inputs of a sampled trajectory under the policy $X\in \mathcal{S}$.  In practice, $T$ is often selected to be finite as the resulted approximation error of \eqref{oracle} decreases exponentially to zero w.r.t. $T$~\cite{malik2019derivative}.  

\subsection{Random Search for Solving (\ref{prim_iterate})} 
{We adopt the random search of Algorithm \ref{alg:learning} to estimate $\nabla_X \mathcal{L}$ via the oracle (\ref{oracle}). The smoothing radius $r$ in Step 4 is used to control its estimation error.} Motivated by~\cite{malik2019derivative}, we shall show that with a large probability, Algorithm \ref{alg:learning} converges and $\{X^{(i)}\}$ remains in the following compact sublevel set
\begin{algorithm}[t]
	\caption{The random search algorithm for (\ref{prim_iterate})}
	\label{alg:learning}
	\begin{algorithmic}[1]
		\Require
		An initial policy $X^{(0)} \in \mathcal{S}$, the number of iterations $N$, a smoothing radius $r$, the stepsize $\eta$, a multiplier $\mu$.
		\For{$i=0,1,\cdots , N-1$}
		\State Sample  $U^{(i)} \in \mathbb{R}^{m\times n}$ uniformly from a unit sphere $\mathbb{S}$ and let $\widehat{X} = X^{(i)}+ rU^{(i)}$.
		\State Obtain a noisy Lagrangian $\widehat{\mathcal{L}}(\widehat{X},\mu)$ from the oracle.
		\State {Compute a gradient estimate} $$\widehat{\nabla_X \mathcal{L}} = \widehat{\mathcal{L}}(\widehat{X},\mu)\frac{n}{r^2}U^{(i)}.$$
		\State Update the policy by $X^{(i+1)} = X^{(i)} - \eta \widehat{\nabla_X \mathcal{L}}$.
		\EndFor	
		\Ensure A policy $X^{(N)}$.
	\end{algorithmic}
\end{algorithm}
%The difficulties in the convergence analysis of Algorithm \ref{alg:learning} hinges on that (a) the gradient dominance property does not hold globally here and (b) the stepsize $\eta$ must be chosen carefully to ensure $X_i \in \mathcal{S}$.
\begin{equation}\label{def:initial set}
\mathcal{S}_{10}=\{X \mid \mathcal{L}(X,\mu)-\mathcal{L}(X^{(0)},\mu) \leq 10 \Delta_{0}\},
\end{equation} 
where $\Delta_0 = \mathcal{L}(X^{(0)},\mu) - \mathcal{L}^*(\mu)$.

Denote 
$
\beta_0 = \sup_{X \in \mathcal{S}_{10}} \beta_X,~~
\xi_0 = \sup_{X \in \mathcal{S}_{10}} \xi_X,~~
\gamma_0 = \inf_{X \in \mathcal{S}_{10}} \gamma_X$, and for notational simplicity, let $\lambda_0=\lambda_{\mathcal{S}_{10}}$(cf. \eqref{dominated}) and $\theta_0 = \min \{1/(2\beta_0),{\gamma_0}/{\xi_0}\}$. Moreover, we make the following assumption in the rest of this section. 

 \begin{assum}\label{assum:uniform}
	The noise sequence $\{w_t\}$ is uniformly bounded, i.e., $\|w_t\| \leq v$, where $v$ is a positive constant.
\end{assum}
Then, we  define 
\begin{equation}\label{def:grad_norm}
\begin{aligned}
&G_{\infty}=\sup _{X \in \mathcal{S}_{10}}\|\widehat{\nabla_X \mathcal{L}}\|_{2},  \\
&G_{2}=\sup _{X \in \mathcal{S}_{10}} \mathbb{E}_{U}[\|\widehat{\nabla_X \mathcal{L}}-\mathbb{E}_{U}[\widehat{\nabla_X \mathcal{L}} | X]\|_{2}^{2}],
\end{aligned}
\end{equation}
both of which are finite under  Assumption \ref{assum:uniform}. 
\begin{theo}\label{theorem:learning}
	Suppose that the stepsize $\eta$ and the smoothing radius $r$ are chosen such that
	\begin{align*}
	\eta &\leq \min \left\{ \frac{\epsilon}{240 \lambda_0 \beta_0 G_2}, \frac{1}{2\beta_0}, \frac{\gamma_0}{G_{\infty}}   \right\} ~~\text{and}\\
	r  &\leq \min\left \{ \frac{\theta_0}{8 \lambda_0 \beta_0 }\sqrt{\frac{\epsilon}{15}}, \frac{1}{2\beta_0} \sqrt{\frac{\epsilon}{30 \lambda_0}}, \gamma_0 \right\}.
	\end{align*}
	For any error tolerance $\epsilon$ such that $\epsilon \log (120 \Delta_{0} / \epsilon)<{10} \Delta_{0}/3$ and $N\ge {4 \lambda_0} \log ({120 \Delta_{0}}/{\epsilon})/{\eta}$, Algorithm \ref{alg:learning} yields that
	\begin{equation}\label{equ:Lag_error}
		\mathcal{L}(X^{(N)}, \mu)- \mathcal{L}^*(\mu) \leq \epsilon
	\end{equation}
with a probability greater than ${3}/{4}$
\end{theo}

The proof is given in Appendix \ref{appd:learn}. In view of \cite{furieri2020learning}, the convergence probability in Theorem \ref{theorem:learning} can be improved to $1-\delta$ for any $0<\delta<1$ by focusing on $
	\mathcal{S}_{\delta}=\{X|\mathcal{L}(X,\mu)-\mathcal{L}(X^{(0)},\mu) \leq 10 \delta^{-1} \Delta_{0}\}.$   

\begin{algorithm}[t]
	\caption{The sample-based primal-dual algorithm for the risk-constrained LQR}
	\label{alg:sample-based}
	\begin{algorithmic}[1]
		\Require  A multiplier $\mu^1 \geq 0$, and a set of stepsizes $\{\zeta^k\}$.
		\For{$k=1,2,\dots$}
		%\State \textbf{Step 1: Learn the dual function}
		\State Solve (\ref{prim_iterate}) by Algorithm \ref{alg:learning} and obtain $\widehat{X}^k$.
		%\State \textbf{Step 2: Dual ascent}
		%\State Compute the subgradient $d^k$ by (\ref{def:sub}).
		%\State Update the multiplier by $\mu^{k+1} = [\mu^{k} + \zeta^k \cdot d^k]_{+}$.
		%\State Learn the policy $X^*(\mu^{k})$ by Algorithm \ref{alg:learning}.
		%\State \textbf{Step 2: Stochastic dual ascent}
		\State Obtain a noisy  $\widehat{J_c}(\widehat{X}^k)$ from the oracle.
		\State Compute a subgradient estimate $\widehat d^k= \widehat{J_c}(\widehat{X}^k) - \bar{\rho}$.
		\State Update the multiplier by $\mu^{k+1} = [\mu^{k} + \zeta^k \cdot \widehat d^k]_{+}$.
		\EndFor	
		%\Ensure Policy-multiplier pair $(X^M,\mu^M)$.
	\end{algorithmic}
\end{algorithm}

\subsection{A Sample-based Primal-dual Algorithm}
In this subsection, we let $\widehat{X}^k=X^{(N)}$ and assume that \eqref{equ:Lag_error} holds  for the sake of simplifying our presentation; see also e.g., \cite{perdomo2021stabilizing}. {The oracle (\ref{oracle}) is adopted to compute a subgradient estimate 
$$
\widehat{d}^k = \widehat{J_c}(\widehat{X}^k) - \bar{\rho}
$$
with the estimation error resulting from the oracle computation (\ref{oracle}) and the gap between $\widehat{X}^k$ and $X^k$. }

{Now, we present our sample-based primal-dual method in Algorithm \ref{alg:sample-based}. Due to the use of biased subgradient estimate, we can obtain the global convergence to a value close to $D^*$.
\begin{theo}\label{theo:mf_pd}
	Let $\bar{\mu}^k = \frac{1}{k} \sum_{i=1}^{k} \mu^{i}$ and $\zeta^k = \text{poly}^{-1}(\epsilon,v)\cdot k^{-1/2}$ in Algorithm \ref{alg:sample-based}, where $\text{poly}(\epsilon,v)$ is a polynomial of degree $4$ and given in Appendix \ref{Appendix:E}. Under Assumptions \ref{assum:slater} and \ref{assum:uniform}, it holds that
	$$
	 D^*-D(\bar{\mu}^k) \leq \text{poly}(\epsilon,v)(k^{-1/2}+4/3).
	$$
\end{theo}}

%\begin{remark}
%	{By Lemma \ref{lem:subgradient}, the convergence rate can be affected by the accuracy $\epsilon$ of primal iteration (\ref{prim_iterate}) and the noise bound $v$.}
%\end{remark}

\section{Simulation}\label{sec:experiment}
In this section, we use simulation to illustrate {the effectiveness of our  RC-LQR}, and the convergence of the policy gradient primal-dual methods in both model-based and sample-based settings.
\subsection{The performance of RC-LQR}
{We adopt the  dynamical model in \cite{tsiamis2020risk}} with
\begin{equation}\label{def:model}
A=\begin{bmatrix}
1 & 0.5 & 0 & 0 \\
0 & 1 & 0 & 0 \\
0 & 0 & 1 & 0.5 \\
0 & 0 & 0 & 1
\end{bmatrix} ~\text{and}~B = \begin{bmatrix}
0.125 & 0 \\
0.5 & 0 \\
0 & 0.125 \\
0 & 0.5
\end{bmatrix}.
\end{equation}
Let
$
Q = \text{diag}(1,0.1,2,0.2)~~\text{and}~~R = I_2.
$
The noise sequence $\{w_t\}$ is given by
$
w_t = \text{clip}(B v_t + e_t),
$
where $v_t=[v_{t,1} ~v_{t,2}]^{\top}$, $e_t$, and clip$(\cdot)$ are chosen as follows.  $\{v_t\}$ is an independent sequence and satisfies that (a) $v_{t,1}$ follows a mixed Gaussian distribution of $\mathcal{N}(3,30)$ and $\mathcal{N}(8,60)$ with weights 0.2 and 0.8, respectively; (b) $v_{t,2}$ follows $\mathcal{N}(0,0.01)$.  $\{e_t\}$ is another Gaussian independent sequence and follows $\mathcal{N}(0,0.01 \times I_4)$. The operator $\text{clip}(\cdot)$ is used to ensure a uniform bound of $w_t$ in Assumption \ref{assum:uniform}, and projects each argurment onto the interval $[-10^4, 10^4]$. Here the statistics of $\{w_t\}$ are evaluated by the Monte Carlo method and the risk tolerance is set as $\bar{\rho} = 15$.
\begin{figure}[t]
	\centerline{\includegraphics[width=60mm]{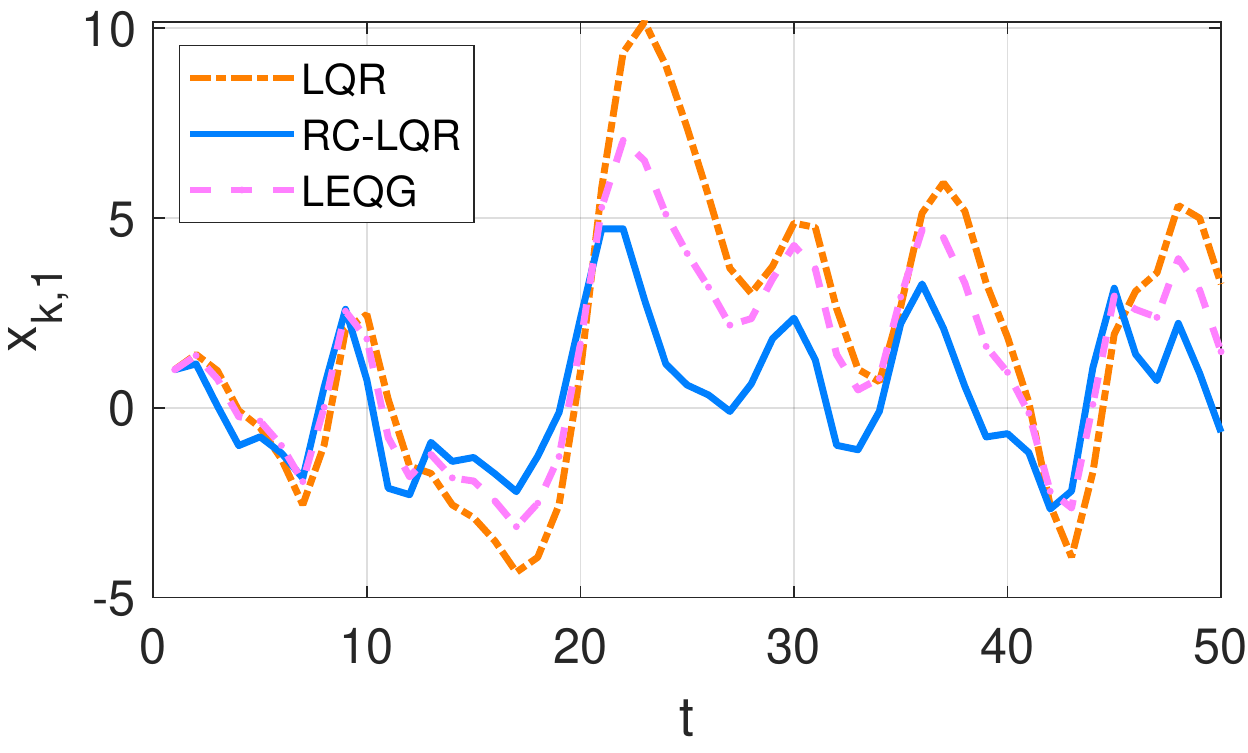}}
	\caption{Evolution of the controlled state $x_{k,1}$ via three different methods.}
	\label{pic:leqg}
\end{figure}

{To illustrate the effectiveness of the RC-LQR, we compare it with the standard LQR and the LEQG with $\theta = 0.01$ \cite{jacobson1973optimal}. Fig. \ref{pic:leqg} depicts the evolution of their controlled states $x_{k,1}$ under the same noise realization, and confirms that our RC-LQR controller compensates the risk better than that of the LQR and LEQG. A similar observation can also be found in \cite{tsiamis2020risk}.}
\begin{figure}[t]
	\centering
	\subfigure[Constant stepsize.]{
		\includegraphics[width=60mm]{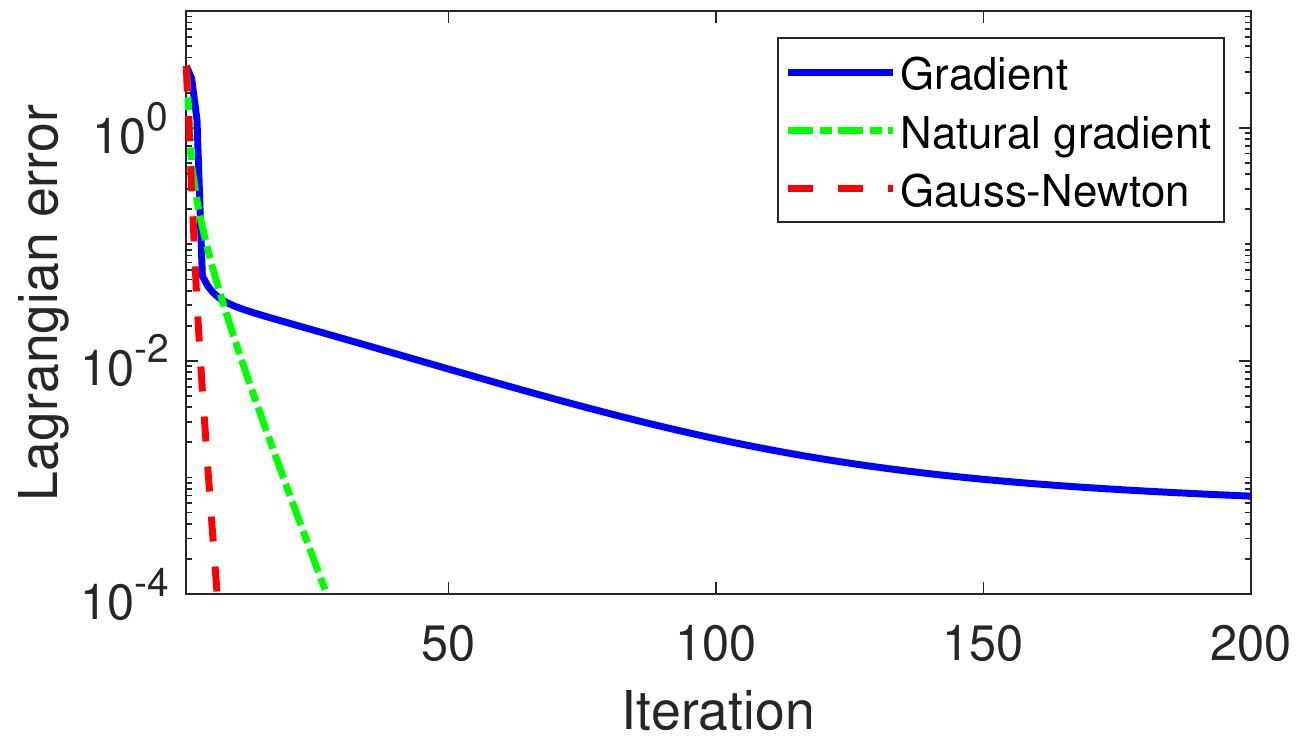}}
	
	\subfigure[Backtracking line search (except GN).]{
		\includegraphics[width=60mm]{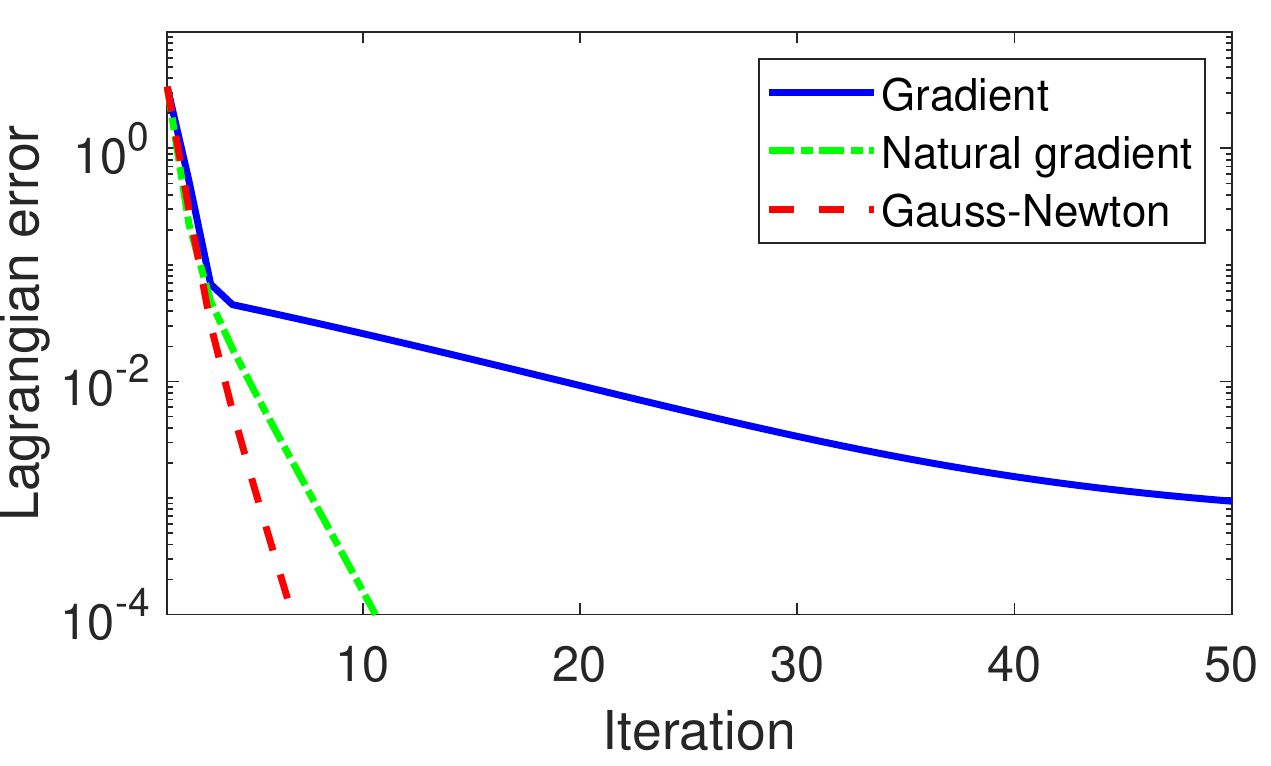}}
	\caption{Relative Lagrangian error  of model-based gradient methods  in (\ref{def:pgmethod}) with $\mu=2$.}
	\label{pic:model_based}
\end{figure}
\subsection{Model-based Setting}
In the model-based setting, we assume that all the parameters in the model (\ref{def:model}) are known. Since the system (\ref{def:model}) is open-loop unstable, we select an initial policy
$$
K^{(0)} = \begin{bmatrix}
0.5 & 0.5 & 0 & 0 \\
0 & 0 & 0.5 & 0.5
\end{bmatrix} ~~\text{and}~~
l^{(0)} = \begin{bmatrix}
-6\\
0
\end{bmatrix}
$$
such that $\rho(A-BK^{(0)})<1$. Since the bounds for the stepsizes in Theorem \ref{theo:pg} could be  conservative in practice, we  manually tune them to be large before  divergence of (\ref{def:pgmethod}) and obtain that $\eta = 3 \times 10^{-3}$ for the PG, $\eta = 0.02$ for the NPG and $\eta = 0.5$ for the GN. We also consider the backtracking line search \cite{boyd2004convex}  with $\alpha = 0.25, \beta = 0.5$ for the PG and NPG where an initial stepsize is set to  $\eta = 0.01$ for the PG and $\eta = 0.05$ for the NPG.  Note that $\eta = 0.5$ is already an optimal stepsize for the GN; see the end of Section \ref{subsec}.

First, we validate the convergence results in Theorem \ref{theo:pg} on the three gradient methods in (\ref{def:pgmethod}).  We adopt the relative Lagrangian error $(\mathcal{L}(X^{(i)},\mu) - \mathcal{L}^*(\mu))/\mathcal{L}^*(\mu)$ to examine the convergence behaviors of (\ref{def:pgmethod}). Fig. \ref{pic:model_based} validates their linear convergence rates of Theorem \ref{theo:pg} and Table \ref{table}.  As expected, it is also observed that the use of a backtracking line search increases the convergence rate.
%\begin{figure}[t]
%	\centering
%	\subfigure[Optimality gap.]{
%		\includegraphics[width=60mm]{optimal_gap_mb}}
%	
%	\subfigure[Risk constraint violation.]{
%		\includegraphics[width=60mm]{risk_mb}}
%	\caption{Convergence of the model-based primal-dual method.}
%	\label{pic:pd_mb}
%\end{figure}

\begin{figure}[t]
	\centerline{\includegraphics[width=60mm]{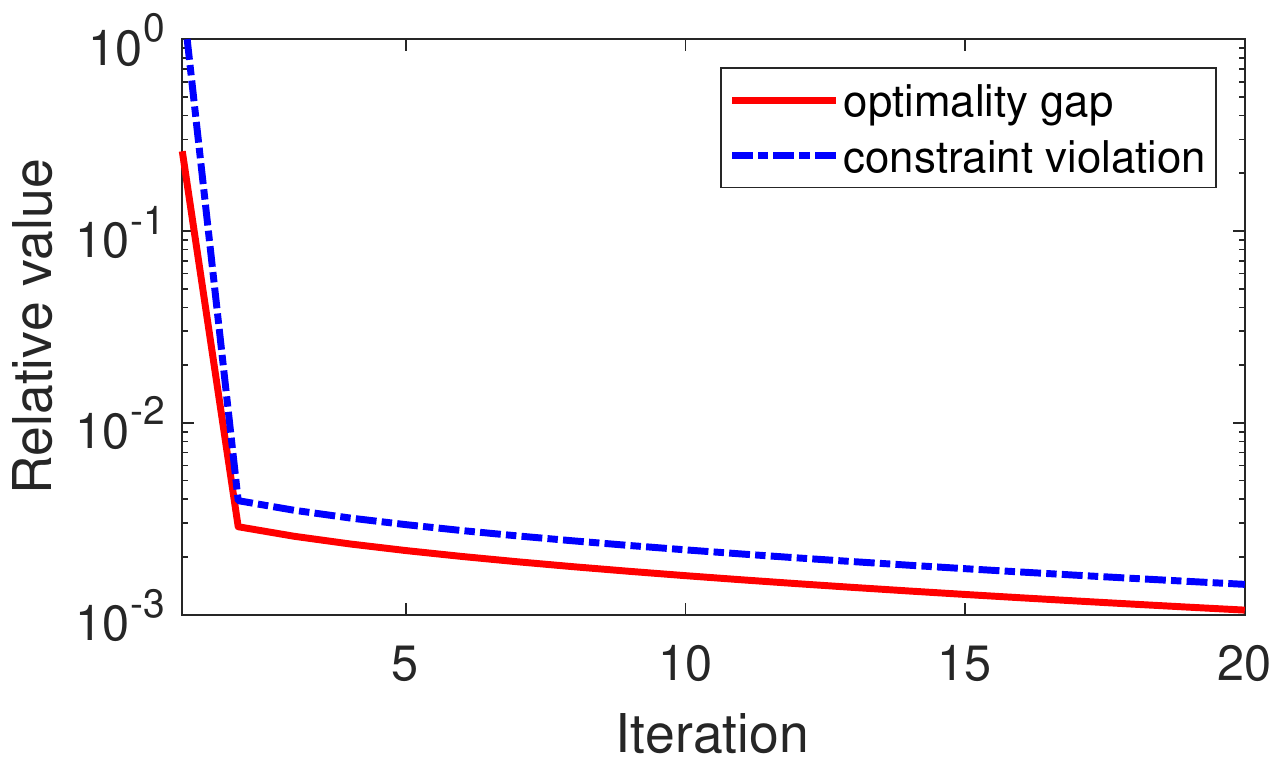}}
	\caption{Convergence of the model-based  Algorithm \ref{alg:model-based}.}
	\label{pic:new}
\end{figure}
Then, we validate the sublinear convergence result in Theorem \ref{theo:model-based}, where the GN is applied to minimize the Lagrangian in Algorithm \ref{alg:model-based}. Let the initial multiplier be $\mu_1 = 0$ and the diminishing stepsize be $\zeta^k = {1}/{(15\sqrt{k})}$. Fig. \ref{pic:new} displays how the relative optimality gap $|J(X^k)- J(X^*)|/J(X^*)$  and the constraint violation $\max\{J_c(X^k)- \bar{\rho},0\}/\bar{\rho}$ decrease to zero. Clearly, both  converge fast under our model-based policy gradient primal-dual method. Note that both the objective function  $J(X)$ and the constraint function $J_c(X)$ are quadratic, and converge with a similar behavior.

%Since by Theorem \ref{theo:model-based} the diminishing stepsize rule leads to a sublinear convergence rate, we additionally perform our primal-dual algorithm with constant stepsize $\zeta^k = 0.1$.

\begin{figure}[t]
	\centerline{\includegraphics[width=60mm]{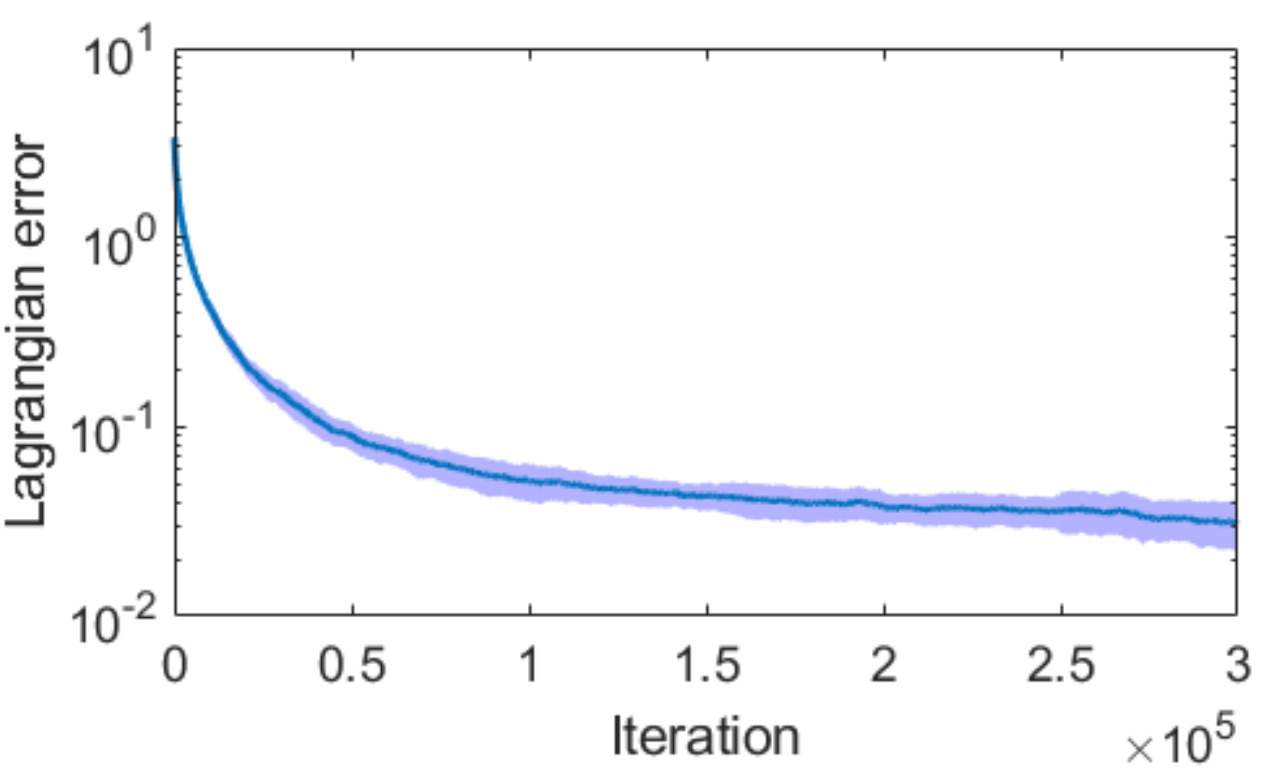}}
	\caption{Relative Lagrangian error of Algorithm \ref{alg:learning} for $\mu = 2$. }
	\label{pic:random_search}
\end{figure}

\begin{figure}[t]
	\centering
	\subfigure[Optimality gap.]{
		\includegraphics[width=60mm]{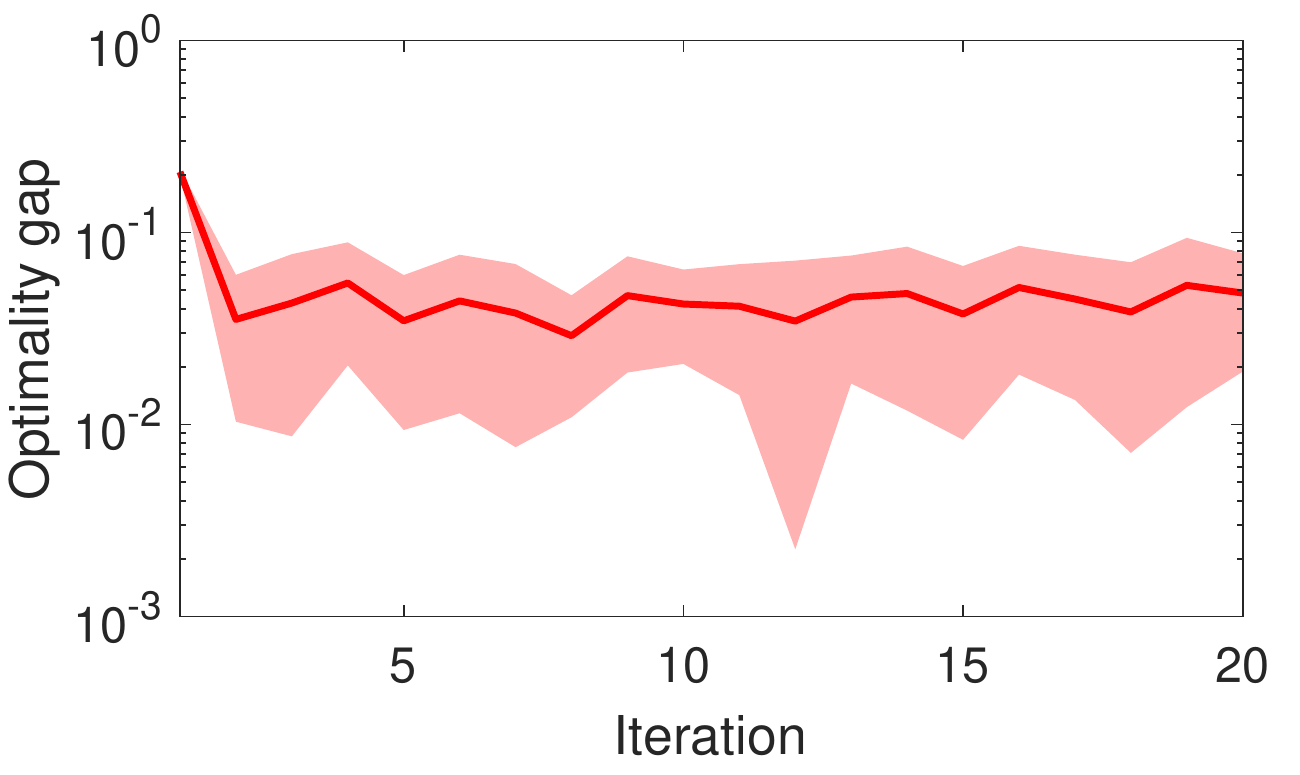}}
	
	\subfigure[Constraint violation.]{
		\includegraphics[width=60mm]{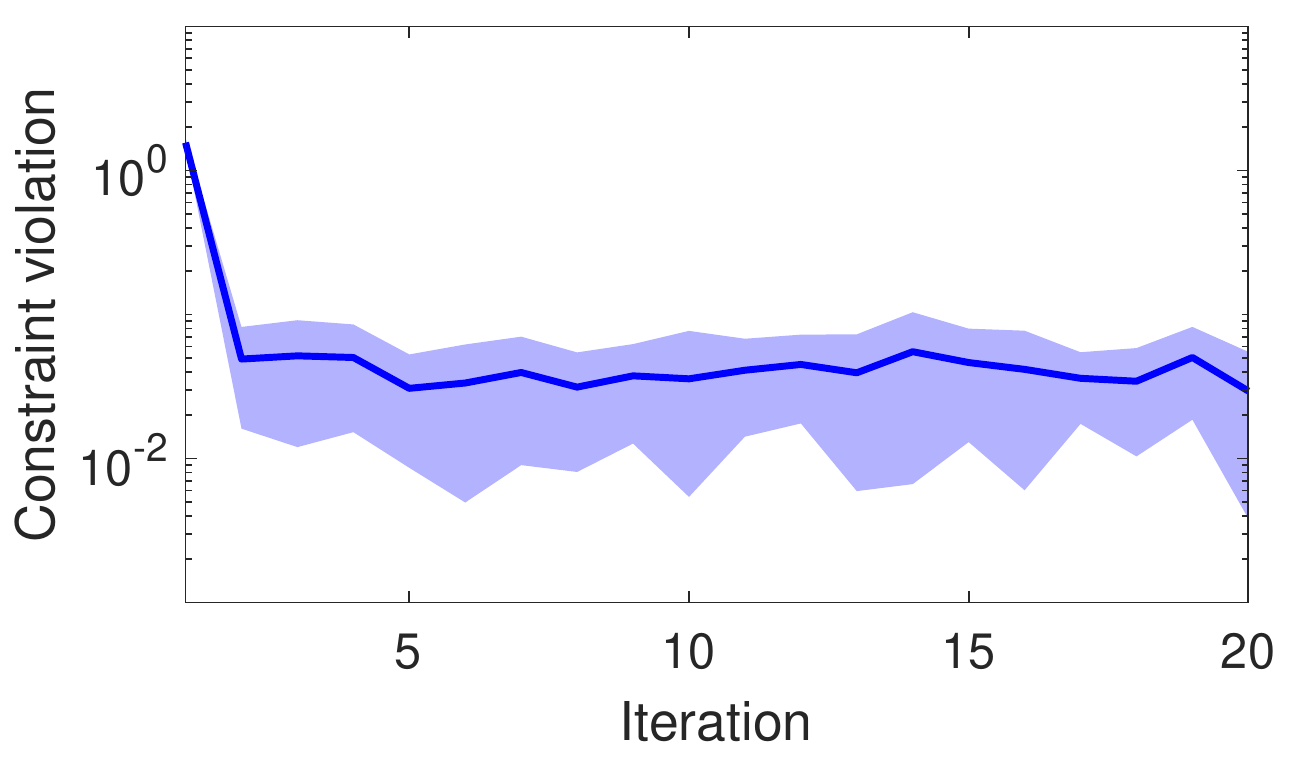}}
	\caption{Convergence of the sample-based Algorithm \ref{alg:sample-based}. }
	\label{pic:pd_mf}
\end{figure}
\subsection{Sample-based Setting}
In the sample-based setting, we use trajectory samples of the system (\ref{equ:sys}) to compute \eqref{oracle} and conduct $20$ independent trials. First, we examine the convergence performance of Algorithm \ref{alg:learning} and set  
the smoothing radius to $r = 0.2$, the sample horizon of the oracle $T=100$ and  the constant stepsize $\eta = 1\times10^{-5}$. Moreover, we display the relative Lagrangian error for $\mu = 2$ in Fig. \ref{pic:random_search}, where the bold centerline denotes the trial mean and the shaded region indicates the variance size. As expected by Theorem \ref{theorem:learning}, Algorithm \ref{alg:learning} converges to a small relative error of $3\%$ with a small variance. 

%\begin{figure}[t]
%	\centerline{\includegraphics[width=60mm]{variance}}
%	\caption{Normalized standard deviation $\textit{std}(\widehat{\mathcal{L}}(X^{(0)},\mu))/\textit{mean}(\widehat{\mathcal{L}}(X^{(0)},\mu))$.}
%	\label{pic:variance}
%\end{figure}

Then, we verify the convergence result of our sample-based primal-dual method in Theorem \ref{theo:mf_pd} by performing Algorithm \ref{alg:sample-based}. Let the initial multiplier be $\mu^1 = 0$ and the diminishing stepsize be $\zeta^k = {1}/{(15\sqrt{k})}$. Fig. \ref{pic:pd_mf} illustrates that both the relative optimality gap and the constraint violation eventually are close to zero.  

%The noises with heavy-tailed distribution may significantly increase the sample complexity of random search in Algorithm \ref{alg:learning}. To illustrate it, we compute $\widehat{\mathcal{L}}(X^{(0)},\mu)$ independently for 100 times, and plot its normalized standard deviation. For comparison, we additionally simulate over a model with small noises $w_{k} \sim \mathcal{N}(0,0.1\cdot I_2)$. Fig. \ref{pic:variance} shows that to achieve the same accuracy, the horizon under heavy-tailed noises is $10^3$ larger than that under small noises.

\section{Concluding Remarks}
In this paper, we have proposed a policy gradient primal-dual framework with global convergence guarantees to solve the RC-LQR problem with a variance-like constraint.  Specifically, we have shown here strong duality, to establish the global convergence, {which in fact can be extended to the case of {\em multiple} constraints}. Such a framework can also be utilized to study linear quadratic tracking.

%We note that our dual iteration has a sublinear convergence rate, which may be sample inefficient in some applications. An important future work is to show the Lipschitz gradient of the dual function, under which the convergence rate can be improved to be linear, based on the perturbation theory for Riccati equations~\cite{sun1998perturbation}. It is also interesting to investigate the sample complexity w.r.t. the number of constraints of our primal-dual policy gradient method.

\section*{Acknowledgement}
We would like to sincerely thank Dr. Kaiqing Zhang from the University of Maryland, College Park, for his constructive comments, and {anonymous reviewers for their valuable suggestions, which significantly helped improve the presentation.}

\bibliographystyle{IEEEtran}
\bibliography{mybibfile}

\appendices

\section{Proofs of some results in Section \ref{sec:primal-dual}}\label{Appendix:A} 

\subsection{Proof of Lemma \ref{prop:lag}}

	Let $u_t = -K x_t + l$. Then, $V_{X}(x)$ satisfies
	{\begin{equation}\notag
	\begin{aligned}
	V_X(x) &= \mathbb{E}\bigg[\sum_{t = 0}^{\infty}\big( x_t^{\top}(Q_{\mu}+K^{\top}RK)x_t  + (2S^{\top}-2l^{\top}RK)x_t \\
	&~~+ l^{\top}Rl - \mu \bar{\rho} -\mathcal{L}(X, \mu)\big) \bigg| x_0 = x \bigg] .
	\end{aligned}
	\end{equation}}
	
	By using backward DP \cite[Chapter 3]{bertsekas1995dynamic}, it can be shown that $V_X(x)$ has a quadratic form~\cite{zhao2021infinitehorizon}, i.e.,
	$
	V_X(x)= x^{\top}P_Kx + g_X^{\top}x + z_X,
	$
	where $P_K,  g_X, z_X$ are to be determined.
	
	{By the Bellman equation~\cite[(3.14)]{sutton1998introduction}, it holds for any $X \in \mathcal{S}$ that
	$
	V_X(x_t) = \mathbb{E}_{w_t}[c_{\mu}(x_{t}, u_{t}) + V_X(x_{t+1})|x_t, u_t=-Kx_t+l]
	$.} That is,
	\begin{align*}
	&x_t^{\top}P_Kx_t + g_X^{\top}x_t + z_X\\
	& = x_t^{\top}(Q_{\mu}+K^{\top}RK)x_t + (2S^{\top}-2l^{\top}RK)x_t + l^{\top}Rl \\
	&+ \mathbb{E} [ ((A-BK)x_t +Bl + w_t)^{\top}P_K ((A-BK)x_t +Bl + w_t) ]\\
	& - \mathcal{L}(X, \mu)- \mu \bar{\rho}+ \mathbb{E}[g_X^{\top}((A-BK)x_t +Bl + w_t)] + z_X\\
	&= x_t^{\top}[Q_{\mu} +K^{\top}RK+(A-B K)^{\top} P_{K}(A-B K)]x_t \\
	&+ 2[-l^{\top}E_K+S^{\top} + \bar{w}^{\top}P_K(A-BK)+g_X^{\top}(A-BK)]x_t \\
	&+ \mathrm{tr}\bigl[P_K(W + (Bl+\bar{w})(Bl+\bar{w})^{\top} )\bigr] \\
	&+ g_X^{\top}(Bl+\bar{w})  + l^{\top}Rl-\mathcal{L}(X, \mu) - \mu \bar{\rho} + z_X.
	\end{align*}
	The proof follows as the equality holds for $x_t \in \mathbb{R}^n$.

\subsection{Proof of Lemma \ref{prop:grad}}
	By Lemma \ref{prop:lag} and observing that $\mathcal{L}(X,\mu)$ is quadratic in $l$, we can compute $\nabla_{l} \mathcal{L}(X,\mu)$ in terms of $G_X$ and $E_K$ as
	\begin{align*}
	\nabla_{l} \mathcal{L}(X,\mu) = 2G_{X} - 2E_K\bar{x}_{X}.
	\end{align*}
	
	We aim to show that
	$
	\nabla_{K} \mathcal{L}(X,\mu) = 2E_K\Sigma_{K} - \nabla_{l}  \mathcal{L}(X,\mu) \bar{x}_{X}^{\top}
	$. First, we express $\mathcal{L}(X,\mu)$ using the stationary distribution $\tau$ of the state as
	\begin{align*}
	&\mathcal{L}(X,\mu)
	= \mathbb{E}_{x \sim \tau} [x^{\top} Q_{\mu} x+ 2x^{\top}S +u_{t}^{\top} R u_{t}- \mu \bar{\rho}] \\
	& =  \mathrm{tr}\{Q_K(\Sigma_{K} + \bar{x}_X\bar{x}_{X}^{\top})\} + 2(S^{\top} -l^{\top}RK)\bar{x}_{X} + l^{\top}Rl- \mu \bar{\rho}.
	\end{align*}
	Then, its gradient in $l$ is given as
	\begin{align*}
	\nabla_l \mathcal{L} & = \nabla_l \operatorname{tr} \{( Q_{\mu}+ K^{\top}RK)\bar{x}_X\bar{x}_X^\top\} \\
	&~~~+ \nabla_{l}\{(2S^{\top} - 2l^{\top}RK)\bar{x}_X\} + 2Rl \\
	& = 2B^{\top}V^{\top}( Q_{\mu}+ K^{\top}RK)\bar{x}_X \\
	&~~~ + 2B^{\top}V^{\top}(S - K^{\top}Rl) -2RK(\bar{x}_X - l).
	\end{align*}

	Similar to \cite[Lemma 1]{fazel2018global}, one can show that
	\begin{align*}
	\nabla_{K} \mathcal{L} &= \nabla_{K} \operatorname{tr}\{(Q_{\mu} + K^{\top}RK)(\Sigma_{K} + \bar{x}_X\bar{x}_X^\top)\}\\
	&~~~ + \nabla_{K}\{(2S^{\top} - 2l^{\top}RK)\bar{x}_X\}\\
	&= 2E_K\Sigma_{K} + 2R(K\bar{x}_X - l)\bar{x}_X^{\top}\\
	&~~~ - 2B^{\top}V^{\top}(( Q_{\mu}+ K^{\top}RK)\bar{x}_X -K^{\top}Rl+S)\bar{x}_X^{\top}\\
	&= 2E_K\Sigma_{K} - \nabla_{l}  \mathcal{L}(X,\mu) \bar{x}_X^{\top}.
	\end{align*}
	Combining
	$
	\nabla_{X} \mathcal{L}(X,\mu) =
	[
	\nabla_{K} \mathcal{L} ~~
	\nabla_{l} \mathcal{L}
	],
	$ and the definition of $\Phi_X$ in (\ref{def:phi}), we have $\nabla_X \mathcal{L}(X,\mu) = 2[
	E_K~~
	G_X
	] \Phi_X$.
	
	Since $\Sigma_{K} > W > 0$, and
	\begin{align*}
	\Phi_X =
	\begin{bmatrix}
	I & - \bar{x}_X \\
	0 & 1
	\end{bmatrix}
	\begin{bmatrix}
	\Sigma_K &  0\\
	-\bar{x}_X^{\top} & 1
	\end{bmatrix},
	\end{align*}
	we have that $\Phi_X$ is strictly positive definite.

\subsection{Proof of Lemma \ref{lem:coer}}
By (\ref{def:bar}), it follows that
$$
\bar{x}_X = (I-A+BK)^{-1}(Bl+\bar{w}).
$$
This implies that $\|\bar{x}_X\| \rightarrow \infty$ if either $\rho(A-BK)\rightarrow 1$ or $\|l\| \rightarrow\infty$, i.e., $\mathcal{L}(X,\mu)$ is coercive over $\mathcal{S}$.

Since the critical point in \eqref{equ:stationary} is unique, we conclude that $\mathcal{S}_{\alpha}$ in (\ref{def:sublevel}) is compact.

\subsection{Proof of Lemma \ref{lem:difference}}\label{appd:difference}
(\romannumeral1) {From the definition of $\mathcal{L}(X,\mu)$ in (\ref{def:L}), it follows that
\begin{align*}
&\hspace{-0.1cm}\mathcal{L}(X',\mu)-\mathcal{L}(X,\mu)=
\lim\limits_{T \rightarrow  \infty} \frac{1}{T} \mathbb{E}\left[ \sum_{t=0}^{T-1}(c_{\mu}(x_t',u_t')-\mathcal{L}(X,\mu))\right]\\
&\hspace{-0.3cm}=\lim\limits_{T \rightarrow  \infty}\frac{1}{T}\left( \mathbb{E}\left[\sum_{t=0}^{T-1}(c_{\mu}(x_t',u_t')-\mathcal{L}(X,\mu))\right] - V_{X}^{T}(x_0')  \right)\\
&\hspace{-0.3cm}= \lim\limits_{T \rightarrow  \infty} \frac{1}{T} \mathbb{E}\left[\sum_{t=0}^{T-1}(c_{\mu}(x_t',u_t')-\mathcal{L}(X,\mu)+ V_{X}^{T}(x_{t+1}')-V_{X}^{T}(x_t'))\right] \\
&\hspace{-0.3cm} = \lim\limits_{T \rightarrow  \infty} \frac{1}{T} \mathbb{E}\left[\sum_{t=0}^{T-1}A_{X}^{T}(x_t', u_t')\right],
\end{align*}
where the expectation is w.r.t. the noise sequence. The second equality follows from the boundedness of $V_{X}^{T}(x_0')$. The third equality follows by telescoping the sum appropriately. The last equality holds by the definition of the advantage function.}

(\romannumeral2) Under $u = -K'x+l'$, we have
\begin{align*}
&\lim\limits_{T \rightarrow  \infty} A_{X}^{T}(x, u)\\
&= c_{\mu}(x,u)-\mathcal{L}(X,\mu) + \mathbb{E}_w[V_{X}(Ax+Bu+w)]-V_{X}(x)\\
& = x^{\top}Q_{\mu}x + 2S^{\top}x + (-K'x+l')^{\top}R(-K'x+l') \\
&-\mathcal{L}(K,l,\mu) + \mathbb{E}_w[V_X((A-BK')x+Bl'+w)] - V_X(x).
\end{align*}
The proof is completed by reorganizing the above terms.
\subsection{Proof of Lemma \ref{lem:gradient dominance}}
It follows from Lemma \ref{lem:difference} that
	\begin{align*}
	&\lim\limits_{T \rightarrow  \infty} A_{X}^{T}(x, u) 
	= ((K'-K)x-l'+l + R_K^{-1}(E_Kx-G_X))^{\top}\\
	&\times R_K ((K'-K)x-(l'-l) + R_K^{-1}(E_Kx-G_X))\\
	&- (E_Kx-G_X)^{\top}R_K^{-1}(E_Kx-G_X) \\
	& \geq - (E_Kx-G_X)^{\top}R_K^{-1}(E_Kx-G_X).
	\end{align*}		
	
	Let $\{x_t^*\}$ and $\{u_t^*\}$ be sequences generated by following the policy $X^*(\mu)$ in (\ref{equ:stationary}). Then, it follows that
	\begin{align}
	&\mathcal{L}(X,\mu)-\mathcal{L}^*(\mu) \notag
	=-
	\lim\limits_{T \rightarrow  \infty} \frac{1}{T}\mathbb{E}\left[\sum_{t=0}^{T-1}A_{X}^{T}(x_t^*, u_t^*)\right]  \notag\\
	& \leq \lim\limits_{T \rightarrow  \infty} \frac{1}{T}\mathbb{E}\left[\sum_{t=0}^{T-1} \operatorname{tr}\{(E_Kx_t^*-G_X)^{\top}R_K^{-1}(E_Kx_t^*-G_X)\}  \right] \notag\\
	& = \lim\limits_{T \rightarrow  \infty} \frac{1}{T}\mathbb{E}\left[\sum_{t=0}^{T-1} \operatorname{tr} \{ 	
	\begin{bmatrix}
	x_t^* \\
	1
	\end{bmatrix}	
	\begin{bmatrix}
	x_t^* \\
	1
	\end{bmatrix}^{\top}	
	[E_K~G_X]^{\top}R_K^{-1}[E_K~G_X]
	\}\right] \notag \\
	& = \operatorname{tr} \{ \Phi^* [E_K~~G_X]^{\top}R_K^{-1}[E_K~~G_X]
	\} \notag\\
	& \leq \|\Phi^*\| \operatorname{tr} \{ [E_K~~G_X]^{\top}R_K^{-1}[E_K~~G_X]
	\} \notag\\
	& \leq \|\Phi^*\| \|R_K^{-1}\|\operatorname{tr} \{ [E_K~~G_X][E_K~~G_X]^{\top} \} \notag\\
	& = ({\|\Phi^*\|}/{\underline\sigma(R)})
	\operatorname{tr} \{ [E_K~~G_X][E_K~~G_X]^{\top} \}. \label{equ:domi}
	\end{align}
	
	 Lemma \ref{prop:grad} implies that $
	 4\operatorname{tr}\{\Phi_X\Phi_X^{\top}[E_K~G_X]^{\top}[E_K~G_X]\}=
	 \operatorname{tr} \{\nabla\mathcal{L}^{\top}\nabla\mathcal{L}\}.
	$ Together with \eqref{equ:domi}, the Lagrangian difference is upper bounded by
	\begin{align*}
	&\mathcal{L}(X,\mu)-\mathcal{L}^*(\mu)\\
	&=\frac{\|\Phi^*\|}{\underline\sigma(R)} \operatorname{tr} \{(\Phi_X \Phi_X^{\top} )^{-1}\Phi_X \Phi_X^{\top}
	[E_K~~G_X]^{\top}[E_K~~G_X]\} \\
	& \leq \frac{\|\Phi^*\|}{4\underline \sigma(R)\underline \sigma(\Phi_X)^2}\operatorname{tr} \{
	\nabla\mathcal{L}^{\top}\nabla\mathcal{L}\}.
	\end{align*}

\subsection{Proof of lemma \ref{lem:cost_diff}}
	Let $\tau'$ be the stationary distribution of the state under $X'$. By Lemma \ref{lem:difference}, it follows that
	\begin{align*}
	&\mathcal{L}(X',\mu)-\mathcal{L}(X,\mu)
	=\lim\limits_{T \rightarrow \infty} \frac{1}{T}\mathbb{E}\left[\sum_{t=0}^{T-1}A_{X}^{T}(x_t', u_t')\right]\\
	&= \mathbb{E}_{x \sim \tau'}\big[((K'-K)x-(l'-l) + R_K^{-1}(E_Kx-G_X))^{\top}\\
	&\times R_K ((K'-K)x-(l'-l) + R_K^{-1}(E_Kx-G_X))\\
	&- (E_Kx-G_X)^{\top}R_K^{-1}(E_Kx-G_X)\big].
	\end{align*}
	
	The rest of the proof follows by (\ref{def:bar}) and (\ref{def:sigma}).

\section{Proof of Lemma \ref{lem:lip}}\label{apx:lip}
We first build a connection between $\mathcal{L}(X, \mu)$ and the standard LQR cost. 
\begin{lemma}\label{lem:cost}
	Define $C(K,\mu) =  \mathrm{tr}(P_KW)$. Then
		\begin{equation}\notag
	C(K,\mu) = \mathrm{tr}\{(Q_{\mu} + K^{\top}RK)\Sigma_{K}\} \leq \mathcal{L}(X,\mu) + S^{\top}Q_{\mu}^{-1}S.
	\end{equation}
\end{lemma}

\begin{proof}
	Comparing the definition of $\mathcal{L}(X,\mu)$ with $C(K,\mu)$, it follows that
	\begin{align*}
		&\mathcal{L}(X,\mu) =\mathrm{tr}\{( Q_{\mu}+ K^{\top}RK)(\Sigma_{K} + \bar{x}_X\bar{x}_X^{\top})\} \\
		& + (2S^{\top} -2l^{\top}RK)\bar{x}_X + l^{\top}Rl- \mu \bar{\rho}\\
		&=\mathrm{tr}\{(Q_{\mu} + K^{\top}RK)\Sigma_{K} \}
		 + (Q_{\mu}\bar{x}_X +S)^{\top}Q_{\mu}^{-1}(Q_{\mu}\bar{x}_X +S)  \\
		& +(K\bar{x}_X -l)^{\top}R(K\bar{x}_X -l)- S^{\top}Q_{\mu}^{-1}S \\
		& \geq C(K,\mu) - S^{\top}Q_{\mu}^{-1}S.
	\end{align*}
\end{proof}
Then, the results in~\cite{fazel2018global} are utilized in our analysis.
Define
\begin{align*}
	&c_1 = \underline{\sigma}(Q_{\mu}) \underline{\sigma}(W)/(4 C(K,\mu)\|B\|(\|A-B K\|+1)), \\
	&\text{and}~c_2  = {(1-\rho(A-BK))}/{(2\|B\|)}.
\end{align*}

\begin{lemma}[\cite{fazel2018global}]\label{lem:fazel}
	Let $X \in \mathcal{S}$. Then, we have that
	
	(\romannumeral1)
	$
	\|\Sigma_{K}\| \leq { C(K,\mu)}/{\underline \sigma(Q_{\mu})}, \|P_K\| \leq {C(K,\mu)}/{\underline \sigma(W)}.
	$
	
	(\romannumeral2)
	If
	$
	\|K^{\prime}-K\| \leq \min (c_1,\|K\|)
	$,
	 it follows that
	\begin{align*}
		&\|P_{K'} - P_K \| \leq  6\|K\|\|R\| ({C(K,\mu)}/({\underline \sigma(Q_{\mu}) \underline \sigma(W)}))^{2}\\
		&\times (\|K\|\|B\|(\|A-B K\|+1)+1)\|K-K'\|.
	\end{align*}
	
	(\romannumeral3) The norm $\|K\|$ can be bounded by
	\begin{equation}\notag
	\|K\| \leq \frac{1}{\underline \sigma(R)}(\sqrt{\frac{\|R_K\|(C(K,\mu)-C(K^{*}(\mu)))}{\underline{\sigma}(W)}}+\|B^{\top} P_{K} A\|).
	\end{equation}
	
	(\romannumeral4)
	$
	{1}/{(1-\rho(A-BK))} \leq {2n\|\Sigma_{K}\|}/{\underline \sigma(W)} \leq {2nC(K,\mu)}/{(\underline \sigma(W) \underline \sigma(Q_{\mu}))}.
	$
\end{lemma}

\begin{lemma}\label{lem:tec_V}
	Let $V' = (I-A+BK')^{-1}$ and
	$
	\|K' - K\| \leq c_2.
	$
	Then,
	$
	\|V' - V\| \leq {2 \|B\| \|K' - K\|}/({1-\rho(A-BK)}).
	$
\end{lemma}

\begin{proof}
	The proof follows from the matrix inverse perturbation theorem~\cite{stewart1990matrix}.
\end{proof}

\begin{lemma} \label{lem:EK}
	For $X \in \mathcal{S}$, we have the following relationship
	$
		\mathrm{tr}(E_K^{\top}E_K) \leq \operatorname{tr} \{ [E_K~~G_X]^{\top}[E_K~~G_X] \} 
		\leq {\|R_K\|}(\mathcal{L}(X,\mu)- \mathcal{L}^*(\mu))/\sigma_0,
	$
	where $\phi_{\alpha}$ is defined in Lemma \ref{lem:local_grad}.
\end{lemma}

\begin{proof}
	Let $X' =X-R_K^{-1} [E_K~~G_X]$, we have
	\begin{equation}\notag
	\begin{aligned}
	&\mathcal{L}(X,\mu)- \mathcal{L}^*(\mu) \geq \mathcal{L}(X,\mu)- \mathcal{L}(X',\mu)\\
	&= \operatorname{tr}\{ [E_K~~G_X]^{\top}R_K^{-1}[E_K~~G_X] \Phi'\}\\
	& \geq {\underline \sigma(\Phi')}\operatorname{tr} \{ [E_K~~G_X]^{\top}[E_K~~G_X] \}/{\|R_K\|} \\
	& \geq {\sigma_0}\operatorname{tr} \{ [E_K~~G_X]^{\top}[E_K~~G_X] \}/{\|R_K\|}.
	\end{aligned}
	\end{equation}
	The last inequality follows since the GN method yields a decrease, i.e., $X' \in \mathcal{S}_0$.
\end{proof}

\begin{lemma}\label{lem:bound_g}
	For all $K'$ such that
	$
	\|K^{\prime}-K\| \leq \min (c_1,c_2),
	$
	we have
	$
	\|	g_{K',l}-g_X\| \leq p_1(K,l) \|K'-K\|
	$, where $p_1(K,l)$ is polynomial in $C(K,\mu)$, $\|A\|$, $\|B\|$, ${\underline{\sigma}(W)}$, ${\underline{\sigma}(Q_{\mu})}$, ${\underline{\sigma}(R)}$,$\|l\|$,$\|\bar{w}\|$.
\end{lemma}

\begin{proof}
	By definition, we have
	\begin{equation}
		\begin{aligned}
		&g_{K',l}^{\top}-g_X^{\top} 
		= 2[-l^{\top}E_{K'}+S^{\top} + \bar{w}^{\top}P_{K'}(A-BK')]V' \\
		& -2[-l^{\top}E_{K}+S^{\top} + \bar{w}^{\top}P_{K}(A-BK)]V\\
		& = -2l^{\top}(E_{K'}-E_K)V'\\
		& + \bar{w}^{\top}(P_{K'}(A-BK') -P_{K}(A-BK))V'\\
		& + 2(-l^{\top}E_{K}+S^{\top} + \bar{w}^{\top}P_{K}(A-BK))(V'-V)\label{gap}.
		\end{aligned}
	\end{equation}
	
	By Lemma \ref{lem:fazel} and the fact
	$
	\|V'\| \leq \|V' - V\| + \|V\|
	\leq{(1+2\|B\| \|K'-K\|)}/{(1-\rho(A-BK))},
	$
	it follows that
	$
	\|2l^{\top}(E_{K'}-E_K)V'\| \leq p(K)\|l\|\|K'-K\|
	$
	where $p(K)$ is polynomial in $C(K,\mu)$, $\|A\|$, $\|B\|$, ${\underline{\sigma}(W)}$, ${\underline{\sigma}(Q_{\mu})}$, ${\underline{\sigma}(R)}$. To bound the second term of \eqref{gap}, we note that
	\begin{align*}
		&\|P_{K'}(A-BK') - P_{K}(A-BK)\|\\
		& \leq \|P_{K'}-P_K\|\|A-BK'\| +\| P_K\|\|B\|K'-K\|.
	\end{align*}
	Again by Lemma \ref{lem:fazel}, it can be bounded by the product of $\|K'-K\|$ and polynomials of related parameters. The third term in \eqref{gap} can be analyzed analogously. Combining the above  completes the proof.
\end{proof}

Next, we derive the Lipschitz constants of $\mathcal{L}(K,l,\mu)$ w.r.t. $K$ and $l$, respectively.
\begin{lemma} \label{lem:K}
	Suppose that
	$
	\|K'-K\| \leq \min (c_1, c_2).
	$
	Then,
	$
	\|\mathcal{L}(K',l,\mu) - \mathcal{L}(K,l,\mu)\| \leq  p_2(K,l) \|K'-K\|,
	$
	where $p_2(K,l)$ is a polynomial in $C(K,\mu)$, $\|A\|$, $\|B\|$, ${\underline{\sigma}(W)}$, ${\underline{\sigma}(Q_{\mu})}$, ${\underline{\sigma}(R)}$,$\|l\|$,$\|\bar{w}\|$, and $\|W\|$.
\end{lemma}

\begin{proof}
	Note that
	\begin{align*}
	&\mathcal{L}(K',l) - \mathcal{L}(K,l)
	=  (g_{K',l}-G_X)^{\top}(Bl+\bar{w})\\
	& +\mathrm{tr}\{(P_{K'}-P_K)(W + (Bl+\bar{w})(Bl+\bar{w})^{\top})\} \\
	&\leq n\|P_{K'}-P_K\|(\|W\| + \|Bl+\bar{w}\|^2)\\ &+\|g_{K',l}-G_X\|\|Bl+\bar{w}\|.
	\end{align*}
	By Lemma \ref{lem:bound_g} and Lemma \ref{lem:fazel}, the proof follows.
\end{proof}

Then, we find the Lipschitz constant of $\mathcal{L}(K,l)$ in $l$.

\begin{lemma} \label{lem:l}
	If $\|l'-l\| \leq \delta$ for some $\delta>0$, then
		$
	\|\mathcal{L}(K,l',\mu) - \mathcal{L}(K,l,\mu)\| \leq  p_3(K,l) \|l'-l\|,
	$
	where $p_3(K,l)$ is a polynomial in $C(K,\mu)$, $\|A\|$, $\|B\|$, ${\underline{\sigma}(W)}$, ${\underline{\sigma}(Q_{\mu})}$, ${\underline{\sigma}(R)}$,$\|l\|$,$\|\bar{w}\|$,$\delta$, $\|W\|$.
\end{lemma}

\begin{proof}
    By direct calculation, we have that
	\begin{equation}\notag
	\begin{aligned}
	&\mathcal{L}(K,l',\mu) - \mathcal{L}(K,l,\mu) \\
	&= (l'-l)^{\top}(R+B^{\top}P_KB)l' + l^{\top}(R+B^{\top}P_KB)(l'-l) \\
	&+2\bar{w}^{\top}P_KB(l'-l) + g_{l'}^{\top}B(l'-l) + (g_{l'}^{\top}-g_l^{\top})(Bl+\bar{w})
	\end{aligned}
	\end{equation}
	and
	$
	\|g_{l'}-g_l\| \leq 2 \| l'-l\|\|E_K\|\|V\|.
	$
	Then, the proof is completed by combining the above lemmas.
\end{proof}

Finally, we prove that the Lagrangian is locally Lipschitz.
\begin{lemma}\label{lem:lag_lip}
	There exist positive scalars $(a_1, b_1)$ that depend on the policy $X = [K~l]$ such that 
	$
	|\mathcal{L}(X',\mu) - \mathcal{L}(X,\mu)| \leq a_1 \| X'-X\|
	$  for all policies $X'$ satisfying $\| X'-X\| \leq b_1$. 
\end{lemma}

\begin{proof}
	By Lemma \ref{lem:K} and Lemma \ref{lem:l}, for $\|X' - X\| \leq \min (c_1, c_2)$ and $\|l'-l\|\leq \delta$, it follows that
	\begin{align*}
		&\mathcal{L}(X',\mu) - \mathcal{L}(X,\mu) \\
		&= \mathcal{L}(K',l',\mu) - \mathcal{L}(K,l',\mu) + \mathcal{L}(K,l',\mu) - \mathcal{L}(K,l,\mu) \\
		&\leq p_2(K,l')\|K'-K\|_F + p_3(K,l)\|l'-l\|_F \\
		&\leq \sqrt{2}\max\{p_2(K,l'), p_3(K,l)\} \|X'-X\|_F,
	\end{align*}
	where $\|\cdot\|_F$ denotes the Frobenius norm.
	
	Letting $a_1 = \sqrt{2}\max\{p_2(K,l'), p_3(K,l)\}$ and $b_1 = \min (c_1, c_2, \delta)$, the proof follows.
\end{proof}

In the sequel, we establish the Lipschitz property for the gradient $\nabla_X \mathcal{L}(X,\mu)$. Similarly, we first derive a bound.
\begin{lemma}\label{lem:bound_phi}
	Suppose that
	$
	\|K'-K\| \leq \min (c_1, c_2).
	$
	Then, it follows that
	$
	\|\Phi_{K',l} - \Phi_X\| \leq p_4(K,l) \|K'-K\|,
	$
	where $p_4(K,l)$ is a polynomial in $C(K,\mu)$, $\|A\|$, $\|B\|$, ${\underline{\sigma}(W)}$, ${\underline{\sigma}(Q_{\mu})}$, ${\underline{\sigma}(R)}$,$\|l\|$,$\|\bar{w}\|$.
\end{lemma}

\begin{proof}
	Note that
	\begin{align*}
	\|\Phi_{K',l} - \Phi_X\| &\leq \operatorname{tr} (\Phi_{K',l} - \Phi_X)\\
	&\leq n\|\Sigma_{K'}-\Sigma_{K}\| +  \|\bar{x}_{K'}\bar{x}_{K'}^{\top} -  \bar{x}_{K}\bar{x}_{K}^{\top}\|.
	\end{align*}
	For the first term, it has been shown in \cite[Lemma 16]{fazel2018global} that
	\begin{equation}\notag
	\|\Sigma_{K^{\prime}}-\Sigma_{K}\| \leq 4(\frac{C(K,\mu)}{\underline{\sigma}(Q_{\mu})})^{2} \frac{\|B\|(\|A-B K\|+1)}{\underline{\sigma}(W)}\|K-K'\|.
	\end{equation}
	Since
$
	\|\bar{x}_{K'}\bar{x}_{K'}^{\top} -  \bar{x}_{K}\bar{x}_{K}^{\top}\| = \|(\bar{x}_{K'}- \bar{x}_{K}) \bar{x}_{K'} +  \bar{x}_{K}(\bar{x}_{K'}- \bar{x}_{K})\|,
$   we have
	\begin{equation}\notag
	\|\bar{x}_{K'}- \bar{x}_{K}\| = \|(V' - V)(Bl+\bar{w}) \| \leq \frac{2 \|B\|\|Bl+\bar{w}\| \|K' - K\|}{1-\rho(A-BK)}
	\end{equation}
	and that
	$
	\|\bar{x}_{K'}\| \leq \|\bar{x}_{K} \| + \| \bar{x}_{K'}- \bar{x}_{K}\| 
	\leq {(2 \|B\| \|K' - K\|+1)\|Bl+\bar{w}\|}/({1-\rho(A-BK)}).
	$
	Then, the proof follows by reorganizing the above terms.
\end{proof}

Next, we establish the Lipschitz constants for $\nabla_{K}\mathcal{L}$ and $\nabla_{l}\mathcal{L}$, respectively.
\begin{lemma}
	Suppose that
	$
	\|K'-K\| \leq \min (c_1, c_2).
	$
	It then follows that
	$
	\|\nabla_{K'} \mathcal{L}- \nabla_{K} \mathcal{L}\| \leq p_5(K,l) \|K'-K\|,
	$
	where $p_5(K,l)$ is a polynomial in $C(K,\mu)$, $\|A\|$, $\|B\|$, ${\underline{\sigma}(W)}$, ${\underline{\sigma}(Q_{\mu})}$, ${\underline{\sigma}(R)}$,$\|l\|$,$\|\bar{w}\|$.
\end{lemma}

\begin{proof}
	We have that
	\begin{align*}
	&\|\nabla_{K'} \mathcal{L}- \nabla_{K} \mathcal{L}\| 
	= 2\|[
	E_{K'}~~
	G_{K',l}
	] \Phi_{K',l} -[
	E_K~~
	G_X
	] \Phi_X\| \\
	&~~~~~~ \leq 2 (\| E_{K'}-E_K\|+\|G_{K',l}-G_X\| )\|\Phi_{K',l}\| \\
	&~~~~~~ + 2(\|[E_K~~G_X]\|)\|\Phi_{K',l}-\Phi_X\|.
	\end{align*}
	
	Note that $\|\Phi_X\|$ can be bounded as
	\begin{equation}\notag
	\begin{aligned}
	\|\Phi_X\| & \leq 1 + \operatorname{tr}(\Sigma_K + \bar{x}_{K}\bar{x}_{K}^{\top}) 
	 \leq 1+ n \|\Sigma_K\| +  \|\bar{x}_{K} \|^2 \\
	& \leq 1 + \frac{n C(K,\mu)}{\underline{\sigma}(Q_{\mu})} + \left(\frac{\|Bl+\bar{w}\|}{1-\rho(A-BK)}\right)^2.
	\end{aligned}
	\end{equation}
	Also, we have
	\begin{align*}
	\|G_{K',l}-G_X\| &\leq \|B\|\|P_{K'} - P_K\|\|Bl + \bar{w}\|\\
	 &+ \|g_{K',l}- G_X\|/2  \leq p_6(K,l) \|K'-K\|,
	\end{align*}
	with $p_6(K,l)$ being polynomial in related parameters. By Lemma \ref{lem:EK}, we obtain that
	\begin{equation}\notag
	\|[E_K~G_X]\| \leq \sqrt{{\|R_K\|}(\mathcal{L}(X,\mu)- \mathcal{L}^*(\mu))/\sigma_0}.
	\end{equation}
	Combining the above inequalities, the proof is completed.
\end{proof}

\begin{lemma}
	For $\|l'-l\|\leq \delta$, we have
	$
	\|\nabla_{l'} \mathcal{L}- \nabla_{l} \mathcal{L}\| \leq p_7(K,l) \|l'-l\|,
	$
	where $p_7(K,l)$ is a polynomial in $C(K,\mu)$, $\|A\|$, $\|B\|$, ${\underline{\sigma}(W)}$, ${\underline{\sigma}(Q_{\mu})}$, ${\underline{\sigma}(R)}$,$\|l\|$,$\|\bar{w}\|$,$\delta$.
\end{lemma}

\begin{proof}
	Note that
	\begin{equation}\notag
	\begin{aligned}
	&\|\nabla_{l'} \mathcal{L}- \nabla_{l} \mathcal{L}\| = 2\|[
	E_{K}~~
	G_{K,l'}
	] \Phi_{K,l'} -[
	E_K~~
	G_X
	] \Phi_X\| \\
	& \leq 2\|G_{K,l'}-G_X\|\|\Phi_{K,l'}\| 
	 + 2\|[E_K~~G_X]\|\|\Phi_{K,l'}-\Phi_X\|.
	\end{aligned}
	\end{equation}
	Then, combining
	$
	\|G_{K,l'}-G_X\| = \| B\|\|g_{K,l'}- G_X\|/2 \leq \| B\|\|E_K\|\|V\|\| l'-l\|,
	$
	and
	\begin{equation}\notag
	\begin{aligned}
	\|\Phi_{K,l'} - \Phi_X\| &\leq  \operatorname{tr}(\Phi_{K,l'} - \Phi_X)
	=  \|\bar{x}_{l'}\bar{x}_{l'}^{\top} -  \bar{x}_{l}\bar{x}_{l}^{\top}\|\\
	& =  	 \|(\bar{x}_{l'}- \bar{x}_{l}) \bar{x}_{l'} +  \bar{x}_{l}(\bar{x}_{l'}- \bar{x}_{l})\| \\
	& \leq (2\|Bl+\bar{w}\|+ \|B\|\delta)\|B\|\|V\|^2\|l'-l\|
	\end{aligned}
	\end{equation}
	completes the proof.
\end{proof}
%Finally, we prove that the gradient is locally Lipschitz.
\begin{lemma}\label{lem:grad_lip}
	There exist positive scalars $(a_2, b_2)$ that depend on the current policy $X = [K~l]$, such that for all policies $X'$ satisfying $\| X'-X\| \leq b_2$, we have
	$
	\|\nabla_{X'}\mathcal{L} - \nabla_{X}\mathcal{L}\| \leq a_2 \| X'-X\|.
	$
\end{lemma}

\begin{proof}
The proof is similar to that of Lemma \ref{lem:lag_lip}.
\end{proof}

So far, we have shown that the Lagrangian and its gradient are locally Lipschitz in Lemma \ref{lem:lag_lip} and Lemma \ref{lem:grad_lip}, respectively. Letting $\gamma_X = \min \{b_1, b_2\}$ and $\xi_{X} = a_1$, $\beta_X = a_2$ completes the proof of Lemma \ref{lem:lip}.

\section{Proof of Theorem \ref{theo:pg}}\label{apx:pg}
\subsection{Proof of the GN update}
	We prove that \text{(\romannumeral1)} under the given stepsize, $X'\in\mathcal{S}$, and \text{(\romannumeral2)} $X'$ stays in the compact sublevel set $S_X = \{ X' |\mathcal{L}(X',\mu) \leq \mathcal{L}(X,\mu) \}$. 
	
	Suppose that \text{(\romannumeral1)} holds (to be proved subsequently). Hence,  $\mathcal{L}(X',\mu)$ is well-defined. By Lemma \ref{lem:cost_diff}, one can show that
	\begin{align*}
	&\mathcal{L}(X',\mu) - \mathcal{L}(X,\mu)\\
	&= (4\eta^2-4\eta) \operatorname{tr} \{ [E_K~G_X]^{\top} R_K^{-1}[E_K~G_X] \Phi'\}\\
	&\leq -2 \eta \operatorname{tr} \{ [E_K~G_X]^{\top}R_K^{-1}[E_K~G_X] \Phi'\}\\
	&\leq -2 \eta \underline\sigma(\Phi') \cdot  \operatorname{tr} \{ [E_K~G_X]^{\top}R_K^{-1}[E_K~G_X] \}\\
	&\leq {-2 \eta \underline\sigma(\Phi')}/{\|\Phi^*\|} ( \mathcal{L}(X,\mu) - \mathcal{L}^*(\mu) ),
	\end{align*}
	where the last inequality follows from (\ref{equ:domi}).
	
	Clearly, it leads to that $\mathcal{L}(X',\mu) \leq \mathcal{L}(X,\mu)$. Thus, $X'$ is contained in $\mathcal{S}_0 \subseteq \mathcal{S}$, and thus
	$$
	\mathcal{L}(X',\mu) - \mathcal{L}(X,\mu) \leq  - ({2 \eta \sigma_0}/{\|\Phi^*\|}) ( \mathcal{L}(X,\mu) - \mathcal{L}^*(\mu) ).
	$$
	
	We have so far shown that for any $0<\eta\leq 1/2$, if the resulting policy $X'$ is stabilizing, then $X' \in S_X$.
	
	To complete the proof, we prove \text{(\romannumeral1)} by contradiction. Suppose that there exists a stepsize $0<\eta'\leq 1/2$ for which the resulting policy is not stabilizing. Consider the ray $\{ X-\eta \nabla_X \mathcal{L} \cdot \Phi_X^{-1}|\eta>0\}$. Let 
	$$\widetilde{\eta} =  \sup \{\eta > 0| \mathcal{L}(X-\eta \nabla_X \mathcal{L} \cdot \Phi_X^{-1}) \leq \mathcal{L}(X, \mu)\}.$$
	Then, it follows from coercivity that there must exist a stabilizing $X'_{\bar{\eta}}=X-\bar{\eta}\nabla_X \mathcal{L} \cdot \Phi_X^{-1}$ in the ray with $\bar{\eta} \in (\widetilde{\eta},\eta')$ such that $\mathcal{L}(X'_{\bar{\eta}})>\mathcal{L}(X,\mu)$. This leads to a contradiction since we can only have $X'_{\bar{\eta}}\in S_X$ by the previous analysis.

\subsection{Proof of the NPG update}
	The stability issue of $X'$ is addressed similarly as in the proof of the GN update. By Lemma \ref{lem:cost_diff}, we obtain that
	\begin{align*}
	&\mathcal{L}(X',\mu) - \mathcal{L}(X,\mu)
	= -4 \eta \operatorname{tr} \{ [E_K~~G_X]^{\top}[E_K~~G_X] \Phi'\}\\
	&+4\eta^{2} \operatorname{tr} \{ [E_K~~G_X]^{\top}R_K[E_K~~G_X] \Phi'\}\\
	&\leq (4\eta^{2} \|R_{K}\| - 4 \eta)\operatorname{tr}  \{ [E_K~~G_X]^{\top}[E_K~~G_X] \Phi'\} \\
	&\leq -2 \eta\operatorname{tr}  \{ [E_K~~G_X]^{\top}[E_K~~G_X] \Phi'\}\\
	&\leq -2 \eta \underline\sigma(\Phi') \operatorname{tr}  \{ [E_K~~G_X]^{\top}[E_K~~G_X] \} \\
	&\leq ({ -2\eta \underline\sigma(\Phi') \underline{\sigma}(R)}/{\|\Phi^*\|}) ( \mathcal{L}(X,\mu) - \mathcal{L}^*(\mu) ),
	\end{align*}
	where the last inequality follows from (\ref{equ:domi}).
	
	Clearly, the Lagrangian decreases as long as $0<\eta\leq {1}/{(2\|R_{K_0}\|)}$. Hence, we obtain that $\underline\sigma(\Phi') \geq \sigma_0$.
	
	Since the iteration
	$
		K' = K - 2\eta E_K
	$
	yields $\|R_{K'}\| \leq \|R_K\|$~\cite{fazel2018global}, it suffices to set the stepsize as $0<\eta\leq {1}/{(2\|R_{K_0}\|)}$. The proof is thus completed.

\subsection{Proof of the PG update}
First, we determine a stepsize $\eta$ such that the Lagrangian decreases after one-step gradient descent.

Define the following quantities
\begin{align*}
	c_3 &=  \frac{\underline{\sigma}(W)\underline{\sigma}(\Phi_X)\underline{\sigma}^2(Q_{\mu})}{48n\|\nabla_K \mathcal{L}\|\|B\|(\|A-B K\|+1)C^2(K,\mu)},\\
	c_4&=\frac{(1-\rho(A-BK))^2\underline{\sigma}(\Phi_X)}{24\|B\|(\|Bl+\bar{w}\|+\|B\|\delta)^2(3-\rho(A-BK))\|\nabla_{K} \mathcal{L}\|},\\
	c_5&=\frac{(1-\rho(A-BK))^2}{12(2\|Bl+\bar{w}\|+ \|B\|\delta)\|B\| \|\nabla_{l} \mathcal{L}\|}, \\
	c_6&= {3}\left({16(1 + n\frac{ C(K,\mu)}{\underline{\sigma}(Q_{\mu})} + (\frac{\|Bl+\bar{w}\|}{1-\rho(A-BK)})^2)\|R_K\| }\right)^{-1}.
\end{align*}

\begin{lemma}\label{lemma:pg}
	Suppose that
	\begin{equation}\label{equ:stepsize}
	\eta \leq \min \left\{\frac{c_1}{\|\nabla_{K}\mathcal{L}\|}, \frac{c_2}{\|\nabla_{K}\mathcal{L}\|}, c_3,c_4,c_5,c_6, \frac{\delta}{\|\nabla_{l}\mathcal{L}\|}\right\}.
	\end{equation}
	Then, the PG update	
	$
	X' =X-\eta \nabla_X \mathcal{L}=X-2 \eta [
	E_K~~
	G_X
	] \Phi_X
	$
	yields that
	$$
	\mathcal{L}(X',\mu) - \mathcal{L}^*(\mu) \leq  (1 -2\eta \frac{  \underline \sigma(\Phi_X)^2 \underline{\sigma}(R)}{\|\Phi^*\|})(\mathcal{L}(X,\mu) - \mathcal{L}^*(\mu)).
	$$
\end{lemma}

\begin{proof}
	The stability issue of $X'$ is addressed similarly as in the proof of the GN update. By Lemma \ref{lem:cost_diff}, one can show
	\begin{equation}\notag
	\begin{aligned}
	& \mathcal{L}(X',\mu)-\mathcal{L}(X,\mu) 
	=-4 \eta \operatorname{tr}(\Phi' \Phi_X [E_K~~G_X]^{\top} [E_K~~G_X])\\
	&+4\eta^{2} \operatorname{tr}( \Phi' \Phi_X [E_K~~G_X]^{\top}R_K [E_K~~G_X]\Phi_X) \\
	&\leq -4 \eta \operatorname{tr}(\Phi_X [E_K~~G_X]^{\top} [E_K~~G_X] \Phi_X)\\
	&+4 \eta\|\Phi'-\Phi_X\| \operatorname{tr}(\Phi_X [E_K~~G_X]^{\top} [E_K~~G_X]) \\
	&+4\eta^{2}\|\Phi'\|\|R_K\| \operatorname{tr}(\Phi_X  [E_K~~G_X]^{\top} [E_K~~G_X]\Phi_X) \\
	&\leq -4 \eta \operatorname{tr}(\Phi_X [E_K~~G_X]^{\top} [E_K~~G_X] \Phi_X)\\
	&+4 \eta \frac{\|\Phi'-\Phi_X\|}{\underline \sigma(\Phi_X)} \operatorname{tr}(\Phi_X [E_K~~G_X]^{\top} [E_K~~G_X] \Phi_X) \\
	&+4\eta^{2} \|\Phi'\|\|R_K\|\operatorname{tr}(\Phi_X [E_K~~G_X]^{\top} [E_K~~G_X] \Phi_X)\\
	&=-4 \eta(1-\frac{\|\Phi'-\Phi_X\|}{\underline \sigma(\Phi_X)}-\eta\|\Phi'\|\|R_K\|) \operatorname{tr}\{\nabla_X \mathcal{L}^{\top} \nabla_X \mathcal{L}\}.
	\end{aligned}
	\end{equation}
	
	By Lemma \ref{lem:gradient dominance}, we obtain that
	\begin{align*}
	&(\mathcal{L}(X',\mu)-\mathcal{L}(X,\mu))/(\mathcal{L}(X,\mu)-\mathcal{L}^*(\mu))\\
	&\leq -4 \eta \frac{\underline \sigma(R)\underline \sigma(\Phi_X)^2}{\|\Phi^*\|}
	\left(1-\frac{\|\Phi'-\Phi_X\|}{\underline \sigma(\Phi_X)}-
	\eta\|\Phi'\|\|R_K\|\right).
	\end{align*}
	
	Thus, it suffices to determine $\eta$ to yield a decreasing cost. To this end, we first bound the norm $\|\Phi'-\Phi_X\|$.
	
	Note that
	\begin{equation}\label{equ:phi}
		\|\Phi'-\Phi_X\| \leq \|\Phi'-\Phi_{K,l'}\| + \|\Phi_{K,l'}-\Phi_X\|.
	\end{equation}

	By Lemma \ref{lem:bound_phi}, if
	$
	\|K'-K\| = \eta \|\nabla_{K} \mathcal{L}\|\leq \min (c_1,c_2),
	$
	then the first term of (\ref{equ:phi}) can be bounded by
	\begin{align*}
	&\|\Phi'-\Phi_{K,l'}\| \\
	&\leq 4n(\frac{C(K,\mu)}{\underline{\sigma}(Q_{\mu})})^{2} \frac{\|B\|(\|A-B K\|+1)}{\underline \sigma(W)}\|K-K^{\prime}\| \\
	&~~~+ \frac{4\|B\|\|Bl'+\bar{w}\|^2}{(1-\rho(A-BK))^2}\|K'-K\|(1+\|B\|\|K'-K\|)\\
	&\leq 4n(\frac{C(K,\mu)}{\underline{\sigma}(Q_{\mu})})^{2} \frac{\|B\|(\|A-B K\|+1)}{\underline \sigma(W)} \|\nabla_{K} \mathcal{L}\| \times \eta\\
	&~~~+ \frac{\|B\|\|Bl'+\bar{w}\|^2}{(1-\rho(A-BK))^2}(6-2\rho(A-BK))\|\nabla_{K} \mathcal{L}\| \times \eta.
	\end{align*}
	
	If $\|l'-l\| = \eta\|\nabla_{l} \mathcal{L}\| \leq \delta$, then the second term of (\ref{equ:phi}) is bounded by
	\begin{equation}\notag
	\begin{aligned}
	\|\Phi_{K,l'} - \Phi_X\|
	 \leq \frac{(2\|Bl+\bar{w}\|+ \|B\|\delta)\|B\| \|\nabla_{l} \mathcal{L}\|}{(1-\rho(A-BK))^2}\times \eta.
	\end{aligned}
	\end{equation}

	Under the given stepsize $\eta$, it can be shown that
	$$
	\|\Phi'-\Phi_X\| \leq ({1}/{12}+{1}/{12}+{1}/{12})\underline \sigma(\Phi_X) =  \underline \sigma(\Phi_X)/4.
	$$
	Then we bound $\|\Phi'\|$ by
	\begin{equation}\notag
	\begin{aligned}
	\|\Phi'\| &\leq \|\Phi'-\Phi_X\| + \|\Phi_X\| \\
	&\leq \frac{1}{4} \underline \sigma(\Phi_X) + 1 + n\frac{ C(K,\mu)}{\underline{\sigma}(Q_{\mu})} + \left(\frac{\|Bl+\bar{w}\|}{1-\rho(A-BK)}\right)^2\\
	&\leq \frac{1}{4} \|\Phi'\| + 1 + n\frac{ C(K,\mu)}{\underline{\sigma}(Q_{\mu})} +\left (\frac{\|Bl+\bar{w}\|}{1-\rho(A-BK)}\right)^2, \\
	\end{aligned}
	\end{equation}
	which implies that 
	\begin{equation}\notag
	\|\Phi'\|  \leq  \frac{4}{3}\left(1 + n\frac{ C(K,\mu)}{\underline{\sigma}(Q_{\mu})} + \left(\frac{\|Bl+\bar{w}\|}{1-\rho(A-BK)}\right)^2\right).
	\end{equation}
	Hence, it follows that
	\begin{equation}\notag
	1-{\|\Phi'-\Phi_X\|}/{\underline \sigma(\Phi_X)}-\eta\|\Phi'\|\|R_K\| \geq 1 - {1}/{4} - {1}/{4} = {1}/{2},
	\end{equation}
	which completes the proof.
\end{proof}

To find a constant stepsize, it suffices to quantify the lower bound of the terms in (\ref{equ:stepsize}). By Lemma \ref{lemma:pg}, we focus on the sublevel set $\mathcal{S}_0 = \{ X |\mathcal{L}(X,\mu) \leq \mathcal{L}(X^{(0)},\mu) \}$. The following inequalities hold:
\begin{align*}
&\underline{\sigma}(\Phi_X) \geq \sigma_0,~ 1-\rho(A-BK) \geq {\underline{\sigma}(W) \underline{\sigma}(Q_{\mu})}/{(2nC(K,\mu))} \\
&{C(K,\mu)}\leq {\mathcal{L}(X,\mu) + S^{\top}Q_{\mu}^{-1}S} \leq {\mathcal{L}(X^{(0)},\mu) + S^{\top}Q_{\mu}^{-1}S}.
\end{align*}		

Since $\mathcal{S}_0$ is compact, $\|l\|$ is bounded by a constant related to $\mathcal{S}_0$. The remaining terms are the gradient norm $\nabla_{l} \mathcal{L}$ and  $\nabla_{K} \mathcal{L}$, which are bounded by $\nabla_{X} \mathcal{L}$.

By the definition of $\nabla_{X} \mathcal{L}$, it follows that
\begin{align*}
\|\nabla_X \mathcal{L}\|^{2} &\leq \operatorname{tr}(\Phi_X [E_{K} ~~G_X]^{\top} [E_{K} ~~G_X] \Phi_X) \\
&\leq\|\Phi_X\|^{2} \operatorname{tr}([E_{K} ~~G_X]^{\top} [E_{K} ~~G_X]).
\end{align*}

By Lemma \ref{lem:EK}, $\|\nabla_X \mathcal{L}\|$ has an upper bound over $\mathcal{S}_0$. Thus, the stepsize $\eta$ has a lower bound polynomial in the problem parameters.

\section{Proof of Theorem \ref{theorem:learning}} \label{appd:learn}

	Our proof is based on~\cite[Theorem 1]{malik2019derivative}. To guarantee the convergence of random search, it requires (a) gradient dominance, (b) locally Lipschitz continuity, and (c) boundedness of gradient norms $G_{\infty} ~\text{and}~ G_{2}$. Thus, we only need to establish (c) by using Assumption \ref{assum:uniform}. Since $G_{2} \leq G_{\infty}^{2}$, it suffices to bound $G_{\infty}$.
	
	We first show that $\widehat{\mathcal{L}}(X, \mu)$ is bounded over $\mathcal{S}_{10}$. By the linear dynamics (\ref{equ:sys}), the state $x_t$ can be written as  
	$$
	x_t = (A-BK)^tx_0 + \sum_{k = 0}^{t-1}(A-BK)^{t-1-k}(Bl+w_t).
	$$
	
By Assumption  \ref{assum:uniform}, it holds that
%	\begin{align*}
%	\widehat{\mathcal{L}}(X, \mu) &= \lim\limits_{T \rightarrow \infty} \frac{1}{T} \sum_{t = 0}^{T}x_t^{\top}(Q_{\mu}+K^{\top}RK)x_t \\
%	&+ 2(S-2K^{\top}Rl)^{\top}x_t + l^{\top}Rl\\
%	&\leq \max \limits_{\{w_t\}} \lim\limits_{t \rightarrow \infty}x_t^{\top}(Q_{\mu}+K^{\top}RK)x_t \\
%	&+ 2(S-2K^{\top}Rl)^{\top}x_t + l^{\top}Rl
%	\end{align*}
	\begin{equation}\label{equ:lagr}
		\begin{aligned}
		\widehat{\mathcal{L}}(X, \mu) 
		&\leq \max \limits_{\{w_t\}} \lim\limits_{t \rightarrow \infty}x_t^{\top}(Q_{\mu}+K^{\top}RK)x_t \\
		&+ 2(S-2K^{\top}Rl)^{\top}x_t + l^{\top}Rl.
		\end{aligned}
	\end{equation}
	
	Noting that $\lim\limits_{t \rightarrow \infty} (A-BK)^t x_0 = 0$, we have that
 	\begin{equation}\label{equ: state_norm}
 		\begin{aligned}
 		\lim\limits_{t \rightarrow \infty}\|x_t\| &\leq \|\sum_{k=0}^{\infty}(A-BK)^k(Bl+w_t)\|\\
 		 &\leq \|Bl+v\|/(1-\rho(A-BK)) \\
 		 & \leq 2nC(K,\mu)\|Bl+v\|/(\underline \sigma(W) \underline \sigma(Q_{\mu})),
 		\end{aligned}
 	\end{equation}
 	where the last inequality follows from Lemma \ref{lem:fazel}.
 	
	Inserting (\ref{equ: state_norm}) into (\ref{equ:lagr}) and noting $\|Q_{\mu}+K^{\top}RK\|\leq \|P_K\|$, then $\widehat{\mathcal{L}}(X, \mu)$ can be bounded by
%	\begin{align*}
%	&\widehat{\mathcal{L}}(X, \mu) \leq \max \limits_{\{w_t\}}(\sum_{k=0}^{\infty}(A-BK)^k(Bl+w_t))^\top\\
%	& \times (Q_\mu + K^{\top}RK)(\sum_{k=0}^{\infty}(A-BK)^k(Bl+w_t)) \\
%	&+ 2(\sum_{k=0}^{\infty}(A-BK)^k(Bl+w_t))^\top (S-K^{\top}Rl) + l^{\top}Rl \\
%	&\leq (\sum_{k=0}^{\infty}\|A-BK\|^k)^2(\|Bl\|+v)^2\|Q_\mu + K^{\top}RK\|\\
%	&+ 2(\sum_{k=0}^{\infty}\|A-BK\|^k)(\|Bl\|+v) \|S-K^{\top}Rl\| + l^{\top}Rl \\
%	&\leq \frac{1}{(1-\rho(A-BK))^2}(\|Bl\|+v)^2\|P_K\|+ l^{\top}Rl\\
%	&+ \frac{2}{1-\rho(A-BK)}(\|Bl\|+v) (\|S\|+\|K\|\|R\|\|l\|),
%	\end{align*}
	\begin{align*}
	&\widehat{\mathcal{L}}(X, \mu) \leq 4n^2C^3(K,\mu)\|Bl+v\|^2/(\underline \sigma^3(W) \underline \sigma^2(Q_{\mu}))+ l^{\top}Rl\\
	&+ 4nC(K,\mu)(\|Bl+v\|)(\|S\|+\|K\|\|R\|\|l\|)/(\underline \sigma(W) \underline \sigma(Q_{\mu})),
	\end{align*} 
	where $C(K,\mu), \|K\|, \|l\|$ are uniformly bounded over $ \mathcal{S}_{10}$. Thus, $\widehat{\mathcal{L}}(X, \mu)$ is bounded. Similarly, it follows from the Lipschitz property that
	$
	\| \widehat{\mathcal{L}}(X+rU,\mu) - \mathcal{L}(X+rU,\mu) \| \leq F
	$
	with some constant $F>0$.
	
	For a given radius $r < \gamma_0$ and a unit perturbation $U \in \mathbb{S}$, the gradient estimate is bounded as
	\begin{align*}
	\|\widehat{\nabla_X \mathcal{L}}\|_{2} &= \frac{n}{r^2}\|\widehat{\mathcal{L}}(X +rU, \mu)\|
	= \frac{n}{r^2}(\|\mathcal{L}(X +rU, \mu)\| + F)\\
	& \leq \frac{n}{r^2} (\|\mathcal{L}(X, \mu)\| + r\xi_0 +F)\\
	& \leq \frac{n}{r^2} (\|10\mathcal{L}(X^{(0)}, \mu)\| + \gamma_0\xi_0 + F).
	\end{align*}
	
	Then, the proof follows directly from \cite[Theorem 1]{malik2019derivative}.

\section{Proof for Theorem \ref{theo:mf_pd}}\label{Appendix:E} 

{Let $\mathcal{S}_{\epsilon}=\{X | \mathcal{L}(X,\mu)-\mathcal{L}^*(\mu) \leq \epsilon\}$.  We first derive a uniform upper bound for $\|\widehat d^k\|$ w.r.t. $\epsilon$ and $v$ via
	\begin{equation}\label{equ:relation}
	\|\widehat d^k\|  \leq \sup_{\{w_t\}, X \in \mathcal{S}_{\epsilon}} \|\widehat{J_c}(X) - \bar{\rho}\|,
	\end{equation}
	which follows from (\ref{equ:Lag_error}). Then, we obtain the following result.
\begin{lemma}\label{lem:subgradient}
	$
	\|\widehat d^k\| \leq p_8(\epsilon,v),
	$ where $p_8(\epsilon,v)$ is a polynomial in $\epsilon$ and $v$ of degree $4$.
\end{lemma}}
	\begin{proof}
	{For any $X \in \mathcal{S}_{\epsilon}$, $\widehat{J_c}(X)$ is bounded by
	\begin{equation}\label{equ:const_bound}
	\begin{aligned}
	&\widehat{J_c}(X) \leq \max \limits_{\{w_t\}}\lim\limits_{t \rightarrow  \infty} (4x_{t}^{\top} QWQ x_{t} + 4x_{t}^{\top}QM_3) \\
	& \leq  \max \limits_{\{w_t\}}\lim\limits_{t \rightarrow  \infty} (4\|x_t\|^2\|QWQ\| + 4\|x_{t}\|\|QM_3\|) \\
	& \leq 4n^2C^2(K,\mu)(\|B\|\|l\|+v)^2\|QWQ\|/(\underline \sigma^2(W) \underline \sigma^2(Q_{\mu})) \\
	& + 8nC(K,\mu)(\|B\|\|l\|+v)\|QM_3\|/(\underline \sigma(W) \underline \sigma(Q_{\mu})),
	\end{aligned}
	\end{equation}
	where the last inequality follows from (\ref{equ: state_norm}).}
{Moreover, it follows from Lemma \ref{lem:cost} that 
	$C(K,\mu) \leq \mathcal{L}(X,\mu) + S^{\top}Q_{\mu}^{-1}S \leq \mathcal{L}^*(\mu) +\epsilon+ S^{\top}Q_{\mu}^{-1}S$. Thus, it suffices to prove an upper bound for $\|l\|$ in (\ref{equ:const_bound}) over $\mathcal{S}_{\epsilon}$.}

{By the definition of $\mathcal{L}(X,\mu)$, it holds that
	\begin{align*}
	&\mathcal{L}(X,\mu) = \mathbb{E}_{x\sim \tau} \left[
	\begin{bmatrix}
	x\\
	-1
	\end{bmatrix}^{\top}
	\begin{bmatrix}
	Q+K^{\top}RK & K^{\top}Rl -S \\
	l^{\top}RK - S & l^{\top}Rl-\mu \bar{\rho}
	\end{bmatrix}
	\begin{bmatrix}
	x \\
	-1
	\end{bmatrix} \right]\\
	& = \text{tr}\{\begin{bmatrix}
	Q+K^{\top}RK & K^{\top}Rl -S \\
	l^{\top}RK - S & l^{\top}Rl-\mu \bar{\rho}
	\end{bmatrix}\Phi_X\} \\
	& \geq \sigma_0\text{tr}\{Q+K^{\top}RK + l^{\top}Rl-\mu \bar{\rho}\} 
	\geq \sigma_0 (\underline{\sigma}(R)\|l\|^2 -\mu \bar{\rho}),
	\end{align*}
	where the first inequality follows from $\mathcal{S}_{\epsilon} \subseteq \mathcal{S}_0$ and the definition $\sigma_0 = \min_{X \in \mathcal{S}_0}  \underline\sigma(\Phi_X)$. Hence, $\|l\|$ is bounded by $
	\|l\|^2 \leq (\mathcal{L}(X,\mu)/\sigma_0 + \mu \bar{\rho})/\underline{\sigma}(R) 
	\leq (\mathcal{L}^*(\mu)+\epsilon + \sigma_0\mu \bar{\rho})/(\sigma_0\underline{\sigma}(R)).$}

{Inserting the bound of $C(K,\mu)$ and $\|l\|$ into (\ref{equ:const_bound}) and noting $J^*(\mu) \leq D^*, \underline{\sigma}(Q_\mu) \geq \underline{\sigma}(Q)$ yields
\begin{align*}
	\widehat{J_c}(X) &\leq \frac{4n^2\|QWQ\|}{\underline \sigma^2(W) \underline \sigma^2(Q)\underline{\sigma}(R)\sigma_0} (D^* + \epsilon + S^{\top}QS)^2\\
	&\times\left({\|B\|}(D^* + \epsilon + \mu \bar{\rho}\sigma_0)^{\frac{1}{2}} + v(\underline{\sigma}(R)\sigma_0)^{\frac{1}{2}} \right)^2 \\
	& + \frac{8n\|QM_3\|}{\underline{\sigma}(W)\underline{\sigma}(Q)(\underline{\sigma}(R)\sigma_0)^{\frac{1}{2}}}(D^* + \epsilon + S^{\top}QS)\\
	&\times \left({\|B\|}(D^* + \epsilon + \mu \bar{\rho}\sigma_0)^{\frac{1}{2}} + v(\underline{\sigma}(R)\sigma_0)^{\frac{1}{2}} \right).
\end{align*}}
{ We note that the bound is also polynomial in $\mu$. As in \cite[Section 4.2]{nedic2009subgradient}, we can focus on a bounded set of $\mu$ since $\mu^*$ in (\ref{def:mu}) is finite, which can be achieved by projection. Thus, without loss of generality, we assume that $\|\mu^k\|\leq e$. Then, $\widehat{J_c}(X)$ is uniformly bounded by a polynomial of $\epsilon$ and $v$, and the proof is completed.}
\end{proof}

Now, we prove Theorem \ref{theo:mf_pd} using standard subgradient arguments; see \cite{nedic2009subgradient}.

By the definition of the projection and subgradient, it holds
\begin{align*}
&\|\mu^{i+1} - \mu^*\|^2 = \| \mu^i - \mu^* + \zeta^i  \widehat{d}^i\|^2 \\
&= \|\mu^i - \mu^*\|^2 +2\zeta^i {d}^i(\mu^i - \mu^*) +2\zeta^i (\widehat{d}^i - d^i)(\mu^i - \mu^*)\\
& + (\zeta^i)^2\|\widehat{d}^i\|^2\\
&\leq \|\mu^i - \mu^*\|^2 +2\zeta^i(D(\mu^i) - D^*) + 8\zeta^i \cdot p_8e + (\zeta^i)^2p_8^2.
\end{align*}
where the inequality follows from Lemma \ref{lem:subgradient} and the boundedness of $\|\mu^k\|$.

Then, rearranging it yields that
$$
D^*-D(\mu^i) \leq \frac{\|\mu^i - \mu^*\|^2}{2\zeta^i} -\frac{\|\mu^{i+1} - \mu^*\|^2}{2\zeta^i} 
+ 4p_8e + \frac{\zeta^ip_8^2}{2}.
$$

Summing up and noting $\zeta^i \geq \zeta^{i+1}$, it follows that
\begin{align*}
&\sum_{i=1}^{k}(D^*-D(\mu^i)) \leq -\frac{1}{2\zeta^{k}}\|\mu^{k+1}-\mu^*\|^2 \\
&+ \frac{1}{{2\zeta^1}}\|\mu^{1} - \mu^*\|^2 +  \frac{1}{2}\sum_{i=1}^{k-1}\left(\frac{1}{\zeta^{i+1}} - \frac{1}{\zeta^i}\right)\|\mu^{i+1}-\mu^*\|^2 \\
&+4p_8e\cdot k+ \frac{p_8^2 }{2}\sum_{i=1}^{k} \zeta^i \leq \frac{2}{\zeta^k}e^2+ 4p_8e\cdot k+\frac{p_8^2 }{2}\sum_{i=1}^{k} \zeta^i.
\end{align*}

By Jenson's inequality, one can easily obtain that
\begin{align*}
D^*-	D(\bar{\mu}^k) \leq \frac{2}{k \zeta^k} e^2+\frac{p_8^2}{2 k} \sum_{i=1}^k \zeta^i + 4be \leq \frac{3p_8e}{\sqrt{k}} + 4p_8e,
\end{align*}
where the last inequality follows by letting $\zeta^k =\frac{1}{p_8e}\sqrt{\frac{2}{k}}$.

\begin{IEEEbiography}[{\includegraphics[width=1in,height=1.25in,clip,keepaspectratio]{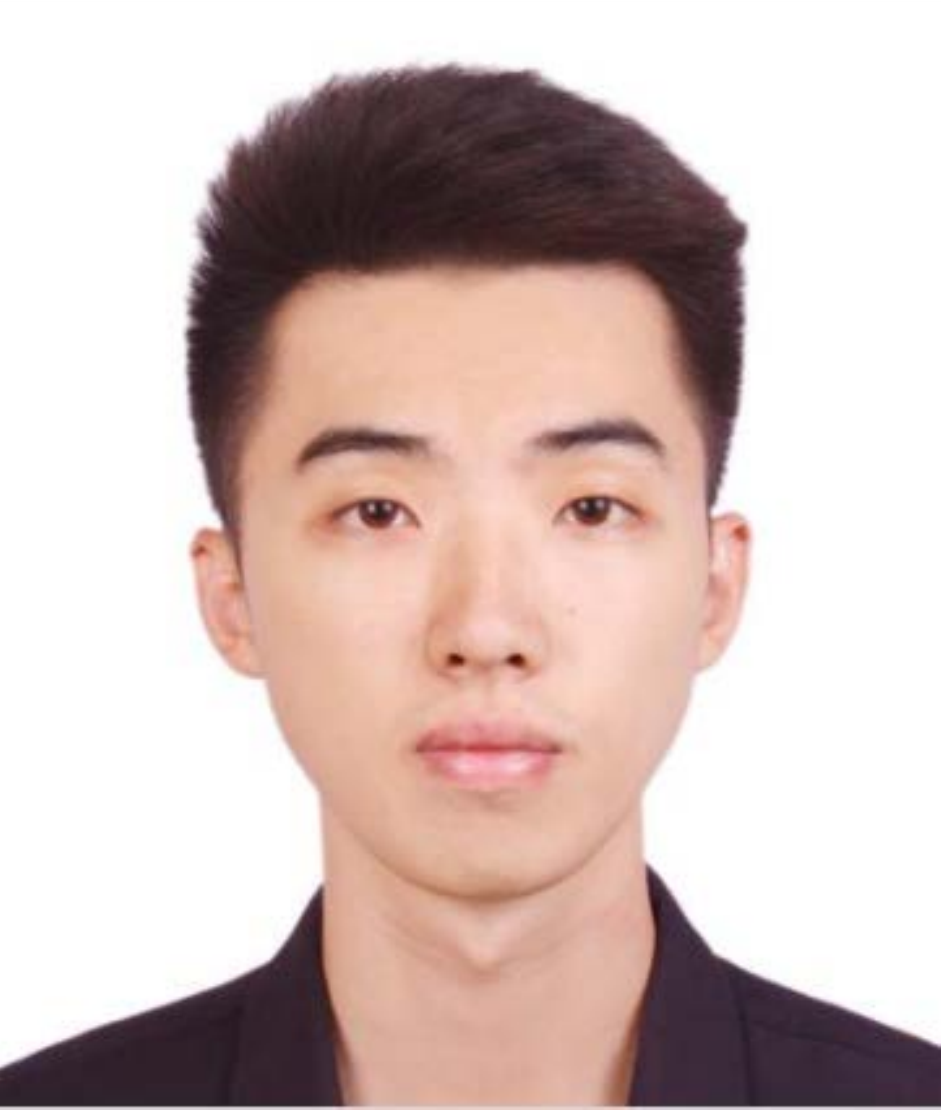}}]{Feiran Zhao} received the B.S. degree in  Control Science and Engineering from the School of Astronautics, Harbin Institute of Technology, Harbin, China, in 2018. He is currently pursuing the Ph.D. degree in  Control Science and Engineering at the Department of Automation, Tsinghua University, Beijing, China. His research interests include reinforcement learning, data-driven methods, control theory and their applications.
\end{IEEEbiography}
\begin{IEEEbiography}
	[{\includegraphics[width=1in,height=1.25in,clip,keepaspectratio]{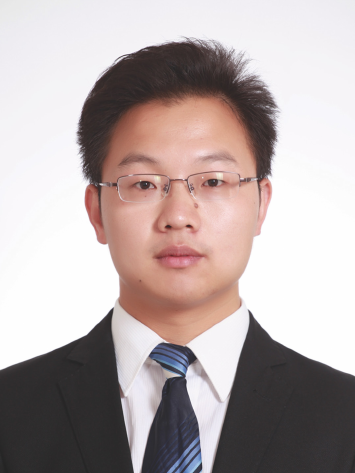}}]
	{Keyou You}(SM'17) received the B.S. degree in Statistical Science from Sun Yat-sen University, Guangzhou, China, in 2007 and the Ph.D. degree in Electrical and Electronic Engineering from Nanyang Technological University (NTU), Singapore, in 2012. After briefly working as a Research Fellow at NTU, he joined Tsinghua University in Beijing, China where he is now a tenured Associate Professor in the Department of Automation. He held visiting positions at Politecnico di Torino,  Hong Kong University of Science and Technology,  University of Melbourne and etc. His current research interests include networked control systems, distributed optimization and learning, and their applications.
	
	Dr. You received the Guan Zhaozhi award at the 29th Chinese Control Conference in 2010,  the ACA (Asian Control Association) Temasek Young Educator Award in 2019 and the first prize of Natural Science Award of the Chinese Association of Automation. He received the National Science Fund for Excellent Young Scholars in 2017. He serves as an Associate Editor for the IEEE Transactions on Control of Network Systems, IEEE Transactions on Cybernetics, and Systems \& Control Letters.
\end{IEEEbiography}

\begin{IEEEbiography}[{\includegraphics[width=1in,height=1.25in,clip,keepaspectratio]{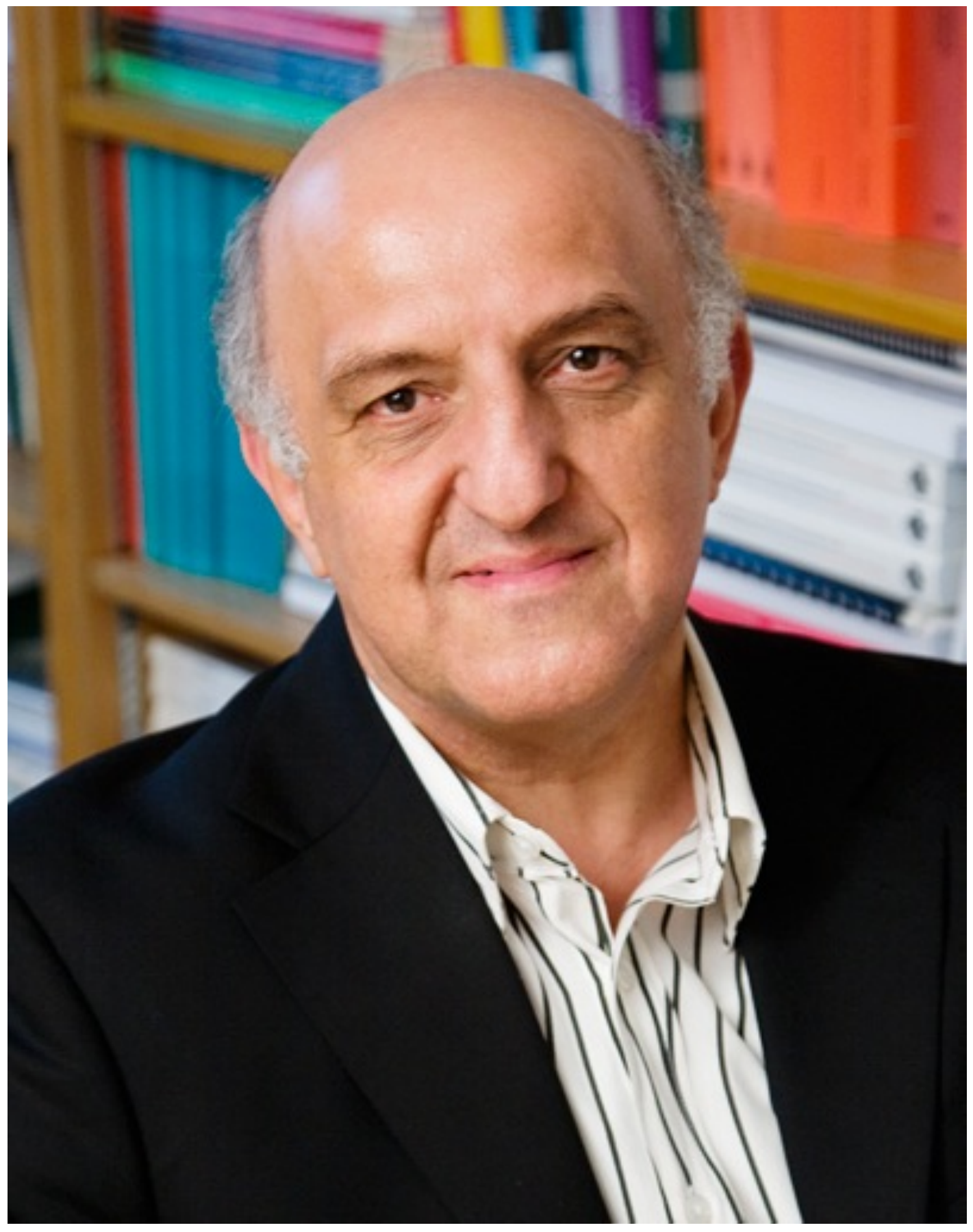}}]
{Tamer Ba\c{s}ar}(S71-M73-SM79-F83-LF13) received the B.S.E.E. degree from the Robert College, \.{I}stanbul, and the M.S., M.Phil., and Ph.D. degrees from Yale University. He has been with the University of Illinois at Urbana-Champaign since 1981, where he is currently Swanlund Endowed Chair Emeritus and Center for Advanced Study (CAS) Professor Emeritus of Electrical and Computer Engineering, with also affiliations with the Coordinated Science Laboratory and the Information Trust Institute. At Illinois, he has also served as Director of CAS (2014-2020), Interim Dean of Engineering (2018), and Interim Director of the Beckman Institute (2008-2010). He has around 1000 publications in systems, control, communications, networks, and dynamic games, including books on non-cooperative dynamic game theory, robust control, network security, wireless and communication networks, and stochastic networked control. His current research interests include stochastic teams, games, and networks; multiagent systems and learning; data-driven distributed optimization; epidemics modeling and control over networks; security and trust; energy systems; and cyber-physical systems.

He is a member of the U.S. National Academy of Engineering and the European Academy of Sciences and Fellow of IEEE, the International Federation of Automatic Control (IFAC), and the Society for Industrial and Applied Mathematics (SIAM). He has received several awards and recognitions over the years, including the highest awards of the IEEE Control Systems Society (CSS), IFAC, the American Automatic Control Council (AACC), and the International Society of Dynamic Games (ISDG), the IEEE Control Systems Award, and a number of international honorary doctorates and professorships. He was the Editor-in-Chief of Automatica from 2004 to 2014. He has served as the President of IEEE CSS, ISDG, and AACC. He is editor of several book series.
\end{IEEEbiography}

\end{document}